\documentclass[11pt]{amsart}

\usepackage{amsmath,amsthm, amscd, amssymb, amsfonts}
\usepackage[all]{xy}
\usepackage[dvips, dvipsnames, usenames]{color}
\usepackage{float}
\usepackage{enumitem}

\allowdisplaybreaks

\usepackage[T2A,T1]{fontenc}
\newcommand\Sh{\mbox{\usefont{T2A}{\rmdefault}{m}{n}\CYRSH}}

\newtheorem{theorem}{Theorem}[section]
\newtheorem*{maintheorem}{Main Theorem}

\newtheorem{lemma}[theorem]{Lemma}
\newtheorem{coro}[theorem]{Corollary}

\theoremstyle{definition}
\newtheorem{definition}[theorem]{Definition}
\newtheorem{example}[theorem]{Example}

\theoremstyle{remark}
\newtheorem{remark}[theorem]{Remark}

\newcommand{\fd}{finite-dimensional}

\newcommand\Irr{\operatorname{Irr}}

\newcommand{\VGamma}{\widehat{\Gamma}}
\newcommand{\VLambda}{\widehat{\Lambda}}

\newcommand{\ydg}{{}^{\ku\Gamma}_{\ku\Gamma}\mathcal{YD}}

\newcommand{\umod}{{}_{\cU}\cM}

\newcommand{\gpl}{\sigma}

\def\G{\mathbb{G}}

\newcommand{\ku}{ \Bbbk}
\newcommand{\kut}{ \Bbbk^{\times}}

\def\ot{\otimes}

\def\Z{\mathbb{Z}}
\def\N{\mathbb{N}}

\def\I{\mathbb I}

\def\If{\mathfrak{I}}

\def\cB{\mathcal{B}}
\def\toba{\mathcal B}

\def\cM{\mathcal{M}}

\def\cC{\mathcal{C}}

\def\cU{\mathcal{U}}
\def\cV{\mathcal{V}}
\def\Uc{\mathcal{U}}

\newcommand{\St}{\mathtt{S}}

\newcommand{\hgamma}{\widehat \Gamma}

\def\ufo{\mathfrak{ufo}}

\def\pf{\begin{proof}}
\def\epf{\end{proof}}

\begin{document}

\title[Simple Modules of the Quantum Double of a Nichols Algebra]{Simple Modules of the Quantum Double of the Nichols Algebra of Unidentified Diagonal Type $\ufo(7)$}

\author[Andruskiewitsch; Angiono; Mej\'\i a; Renz]
{Nicol\'as Andruskiewitsch, Iv\'an Angiono,  Adriana Mej\'\i a, Carolina Renz}

\address{FaMAF-CIEM (CONICET), Universidad Nacional de C\'ordoba,
	Medina A\-llen\-de s/n, Ciudad Universitaria, 5000 C\' ordoba, Rep\'
	ublica Argentina.} \email{(andrus|angiono)@famaf.unc.edu.ar}

\address{CAPES-PNPD, Departamento de matem\'aticas / Campus Trindade. Universidade Federal de Santa Catarina, Florianopolis, SC, Brasil 88040-900} \email{sighana25@gmail.com}

\address{Universidade do Vale do Rio dos Sinos, Av. Unisinos, 950, Cristo Rei, 93022-000, São Leopoldo - RS, Brasil} \email{carolinacnr@gmail.com}

\thanks{The work of N. A., I. A. and A. M. was partially supported by CONICET,
	Secyt (UNC), the MathAmSud project GR2HOPF. The work of I. A. was partially supported by ANPCyT (Foncyt).}

\begin{abstract}
The \fd{}  simple modules over the Drinfeld double of \\ the bosonization of the Nichols algebra $\ufo(7)$
are classified.
\end{abstract}

\maketitle
%\tableofcontents
\section*{Introduction}

\subsection{Motivations and context} 
The purpose of this paper is to compute explicitly all simple \fd{} modules
of the  Hopf algebra $\cU$ that we introduce explicitly below
by generators and relations  after some necessary notation.

In short, $\cU\simeq D(H)$ arises as the Drinfeld double of $H =\toba(V)\#\ku\Lambda$, where 
$\Lambda$ is an abelian group, $V$ is a specific braided vector space of diagonal type 
(realized as a Yetter-Drinfeld module over $\Lambda$)
and $\toba(V)$ denotes its Nichols algebra. 
The representation theory of such Drinfeld doubles or slight variations thereof was treated in many papers,
among them \cite{L,AJS,RS,ARS,H-lusztig iso,HY,PV}. 
Indeed, the first two articles deal with the representation theory of the finite quantum 
groups or Frobenius-Lusztig kernels, while in the others some general results are established.
Presently we know that simple modules are parametrized by highest weights but we ignore the
character formulas and  the dimensions, in the general situation. 
The goal of working out this particular toy example, establishing the dimensions of all
irreducible modules of $\Uc$, is to gain experience for further developments. 
The algebra $\cU$ that we are dealing with is small enough to allow an approach by elementary computations.
Arguing as in \cite[Theorem 3.7]{ARS}, see also \cite[Proposition 5.6]{H-lusztig iso}, 
it is possible to prove that $\Uc$ is a quasi-triangular Hopf algebra, even a ribbon one by
the criterion in \cite[Theorem 3]{RK}, what makes it susceptible of applications.
If $\Lambda$ is finite, then the simple $\cU$-modules are just the simple Yetter-Drinfeld $H$-modules;
therefore the classification here might have applications to the study of basic Hopf algebras.

Why do we choose this specific $V$? Let us first recall that 
finite-dimensional Nichols algebras of diagonal type were classified in \cite{H-classif}.
We find it useful to organize the classification in four classes:
\begin{itemize}
	\item Standard type \cite{A-standard}, including Cartan type \cite{AS Fin-qg}.
	\item Super type \cite{AAY}.
	\item (Super) modular type \cite{AA}.
	\item Unidentified type \cite{A-unidentified}.
\end{itemize}
Actually, $\toba(V)$ is the smallest Nichols algebra of unidentified type. 
Following the terminology introduced in \cite{AA},
we say that $V$ is of type $\ufo(7)$.
Indeed, the Nichols algebra $\toba(V)$ is \fd{}
by \cite[Table 1, row 7]{H-classif}; more precisely, cf. \eqref{def:PBW1},
\begin{align*}
\dim\cB(V)=2^43^2=144.
\end{align*}
By \cite{A-unidentified},  a consequence of \cite{A-convex, A-presentation}, we know that
$\cB(V)$ has a presentation by generators $E_1$, $E_2$ and relations \eqref{th:relations double nichols-E} below.
Thus $\cB(V)$ is manageable yet does not arise from any Lie algebra, what makes it attractive.

\medbreak
There is another reason to address the representation theory of $\cU$.
A \fd{} Nichols algebra of diagonal type admits both a distinguished pre-Nichols algebra \cite{A-distinguished} and a distinguished post-Nichols algebra \cite{AAR}; the representation theories
of the corresponding Drinfeld doubles seem to be very rich. However our $\toba(V)$ coincides
with its distinguished pre-Nichols and post-Nichols algebras, being therefore of singular interest 
(the only  other Nichols algebra with this feature has diagram 
$\xymatrix{{\circ}^{\omega} \ar@{-}[r]^{-\omega} & {\circ}^{-1} }$, $\omega\in \G_3'$, which is of standard type $B_2$).
This peculiar behaviour appeals to the consideration of $V$. 

\subsection{The algebra $\Uc$}
We now introduce formally $\Uc$. Let us begin with some notation.

 If $k, \ell \in \N_0$, then we denote $\I_{k, \ell} = \{n\in \N_0: k\le n \le \ell \}$; also $\I_{\ell} := \I_{1, \ell}$.
Let $\ku$ be an algebraically closed field  of characteristic zero and $\ku^{\times} = \ku - 0$.
Let $\G_{12}$  be
the group of $12$-roots of unity in $\ku$, and let  $\G'_{12}$ be the subset of primitive roots of order $12$.

\medbreak
To define $\Uc$, we need some data:

\begin{itemize}[leftmargin=*]\renewcommand{\labelitemi}{$\circ$}
	\item A matrix $\mathbf{q} = (q_{ij})_{1\le i,j\le 2} = \begin{pmatrix} \zeta^4 & q_{12}  \\ q_{21} & -1   \end{pmatrix} \in \ku^{2\times 2}$ such that $q_{12}q_{21} = \zeta^{11}$; that is,
	 its associated generalized Dynkin diagram is given by
		\begin{align}\label{eq:diagram-used}
		\xymatrix{\underset{1}{\circ}^{\zeta^4} \ar@{-}[r]^{\zeta^{11}} &
			\underset{2}{\circ}^{-1} }.
		\end{align}

\medbreak \item	An abelian group $\Lambda$ whose group of characters
is denoted by $\VLambda$. We set $\Gamma = \Lambda \times \VLambda$.

\medbreak \item  $g_1$, $g_2 \in \Lambda$, $\gpl_1$, $\gpl_2 \in \widehat{\Lambda}$
	such that $ \begin{pmatrix} \gpl_1(g_1) & \gpl_2(g_1)  \\ \gpl_1(g_2) & \gpl_2(g_2)   \end{pmatrix} =
	\begin{pmatrix} \zeta^4 & q_{12}  \\ q_{21} & -1   \end{pmatrix}$.
	\end{itemize}

\medbreak
Starting from this data, we consider vector spaces $V$ and $W$ with bases $v_i$, respectively $w_i$, $i\in \I_2$ and define an action and a 
$\Gamma$-grading on $V$ and $W$ by
\begin{align}\label{eq:intro-actions}
g \cdot  v_i&= \gpl_i(g) v_i, & \gpl \cdot v_i&= \gpl(g_i) v_i,&
g\cdot  w_i&= \gpl_i^{-1}(g) w_i, & \gpl \cdot  w_i&= \gpl(g_i^{-1}) w_i;
\\ \label{eq:intro-coactions}
\deg  v_i&= g_i, &
\deg  w_i&= \gpl_i, & g&\in \Lambda, \, \gpl\in \VLambda,  & i&\in \I_2.
\end{align}
Then $V\oplus W$ is a Yetter-Drinfeld module over $\ku \Gamma$
and $T(V\oplus W)$ is a braided Hopf algebra in $\ydg$.
In particular $V$ is a braided vector space of diagonal type $\ufo(7)$, as said.

\medbreak
It is convenient to start with the auxiliary Hopf algebra $\overline{\cU} = T(V\oplus W) \# \ku \Gamma$;
in particular, $T(V\oplus W)$ and $\ku \Gamma$ are subalgebras of $\overline{\cU}$ and
\begin{align*}
g  v_i&= \gpl_i(g) v_i, & \gpl  v_i&= \gpl(g_i) v_i \gpl,&
g  w_i&= \gpl_i^{-1}(g) w_ig, & \gpl   w_i&= \gpl(g_i^{-1}) w_i\gpl,
\end{align*}
$g\in \Lambda$, $\gpl\in \VLambda$,  $i \in \I_2$. To stress the similarity with quantum groups,
 we denote in $\overline{\Uc}$ or any quotient thereof, as in \cite{ARS,H-lusztig iso,HY},
\begin{align}
E_i&=v_i,& F_i&=w_i\sigma_i^{-1}, & i&\in \I_2.
\end{align}
Thus 	
\begin{align*}
g E_i&= \gpl_i(g) E_i g,& \gpl E_i&= \gpl(g_i) E_i\gpl,& g F_i &= \gpl_i^{-1}(g)F_i g,&
\gpl F_i &= \gpl(g_i^{-1}) F_i\gpl.
\end{align*}

We also need the notation of the so-called root vectors, needed for the relations and for the PBW-basis:
\begin{align*}
E_{12}&= E_1E_2 - q_{12}E_2E_1,& E_{112} &=E_1E_{12} - q_{12}\zeta^4E_{12}E_1,  &
E_{11212} &= E_{112}E_{12} - q_{12}\zeta E_{12}E_{112},\\
F_{12}&= F_1F_2 - q_{21}F_2F_1,& F_{112} &=F_1F_{12} - q_{21}\zeta^4F_{12}F_1,  &
F_{11212} &= F_{112}F_{12} - q_{21}\zeta F_{12}F_{112}.
\end{align*}

We are now ready to define $\Uc$.

\begin{definition} The  algebra $\Uc$ is the quotient of $\overline{\Uc}$ by the ideal generated by
	\begin{align}\label{th:relations double nichols-E}
	E_1^2&=0,& E_2^2&=0,&  E_{11212}E_{12} &= \zeta^{10}q_{12}E_{12}E_{11212}, \\
\label{th:relations double nichols-F}	F_1^2&=0,& F_2^2&=0,&  F_{11212}F_{12} &= \zeta^4q_{21}F_{12}F_{11212},
\\ \label{th:relations double nichols-EF}
&&	E_kF_i-F_iE_k &= \delta_{ki}(g_i-\sigma_i^{-1}).&&
	\end{align}
The algebra $\cU$ is a Hopf algebra with coproduct given by
	\begin{align*}
	\Delta(E_i)&=E_i\ot 1+ g_i\ot E_i,& \Delta(F_i)&=F_i\ot \gpl_i^{-1} +1\ot F_i,& \Delta (g)&=g\ot g,& g\in \Gamma.
	\end{align*}
\end{definition}

Let $\Uc^-$ (respectively $\Uc^+$) be the subalgebra of $\Uc$ generated by $F_1, F_2$
(respectively $E_1, E_2$). The following facts are not difficult to prove and can be derived from 
general results in the literature cited above:

\begin{itemize}[leftmargin=*]\renewcommand{\labelitemi}{$\circ$}
	\smallbreak \item $\Uc$ has a triangular decomposition
	$\Uc\simeq \Uc^+\otimes\ku \Gamma\otimes \Uc^-$, given by the multiplication map.
	
	\smallbreak \item $\Uc^+ \simeq \toba(V)$; in what follows we identify these two algebras.
	
	\smallbreak \item $\cU$, $\Uc^+$ and $\Uc^-$,  admit a $\Z^{2}$-graduation $\cU = \oplus_{\beta \in \Z^2} \cU_{\beta}$ such that $\deg E_i=\alpha_i = -\deg F_i$, $i\in \I_2$, and  $\deg x =0$ for $x\in \Gamma$.	
\end{itemize}

Here $(\alpha_i)_{i\in \I_2}$ is the canonical basis of $\Z^{2}$.

\subsection{Verma modules}
We recall succinctly the description of the simple modules in terms of highest weights.

Let $\umod$ be the category of  left $\cU$-modules
and let $\Irr \Uc$ be the set of isomorphism classes of \fd{} simple $\Uc$-modules.
If $M \in \umod$ and $\lambda \in \hgamma$, then 
\begin{align*}
M^{\lambda} = \{m\in M: g\cdot m = \lambda(g)m \ \forall g \in \Gamma \}
\end{align*}
is the space of weight vectors with weight $\lambda$; if $M = \oplus_{\lambda\in \hgamma} M^{\lambda}$, then we say that $M$ is diagonalizable.

Let $\lambda\in \widehat{\Gamma}$. We denote by $\ku_{\lambda}$ the  $\ku \Gamma\otimes\Uc^-$-module
defined by $\lambda \otimes \varepsilon$ (the counit).
The \emph{Verma module} $M(\lambda)$ associated to $\lambda$ is the induced module
\begin{align}
M(\lambda)=\text{Ind}_{\ku \Gamma\otimes \Uc^-}^{\Uc} \ku_{\lambda} \simeq \Uc/\big(\Uc F_1+\Uc F_2+\sum_{g\in\Gamma}\Uc(g-\lambda(g))\big).
\end{align}
Let $v_\lambda$ be the residue class of $1$ in $M(\lambda)$; then 
we have an isomorphism of $\Uc^+$-modules
\begin{align*}
\Uc^+ &\simeq M(\lambda), &1 &\longmapsto v_\lambda.
\end{align*}
Hence $\dim M(\lambda) = \dim \toba(V) = 144$.
Thus the PBW-basis of $\Uc^+ \simeq \toba(V)$ becomes via this isomorphism 
a basis of $M(\lambda)$.

The $\Z^2$-grading on $\Uc^+ \simeq \toba(V)$ induces a $\Z^2$-grading on $M(\lambda)$ such that
\begin{align*}
M(\lambda)_\beta&=\Uc_\beta\cdot v_\lambda, & \beta&\in\Z^2.
\end{align*}
Thus 
\begin{align*}
M(\lambda)_0 &=\ku v_\lambda,& 
\Uc_\beta\cdot M(\lambda)_\gamma &\subset M(\lambda)_{\beta+\gamma},&  \beta,\gamma&\in\Z^2.
\end{align*}
The family of $\Uc$-submodules of $M(\lambda)$ contained in $\sum_{\beta\neq 0}M(\lambda)_\beta$ has a unique maximal element $N(\lambda)$. We set
\begin{align*}
 L(\lambda)= M(\lambda)/N(\lambda).
\end{align*}

Since $\Uc$ satisfies the conditions on \cite[Section 2]{RS},  \cite[Corollary 2.6]{RS} implies that
\begin{align}\label{eq:simple-mod}
&\text{\emph{The map} }& \lambda&\mapsto L(\lambda)& &\text{\emph{ provides  a bijective correspondence}  } & \hgamma &\simeq \Irr \Uc.
\end{align}

Alternatively  we see that $L(\lambda)$ is simple arguing as in \cite[Theorem 1]{PV}; then \cite[Theorem 3]{PV} gives \eqref{eq:simple-mod}.
Notice that $L(\lambda)$ inherits the grading from $M(\lambda)$. Also, it follows that every simple $M\in \umod$ is diagonalizable. 

\medbreak
Lowest weight modules of weight $\lambda$ are defined as usual; $M(\lambda)$ covers every lowest weight module of weight $\lambda$, that
in turn covers $L(\lambda)$. Highest weight modules are defined similarly.

\subsection{Main result}
In our main theorem,  we give the dimension of $L(\lambda)$ for each $\lambda\in \widehat{\Gamma}$,
in terms of certain equalities arising from the Shapovalov determinant \cite{HY} satisfied by
\begin{align*}
\lambda_i&=\lambda(g_i \gpl_i),& i&\in \I_2.
\end{align*}
Indeed the Shapovalov determinant in the context of this paper is
\begin{align} \label{eq:shapovalov}
\begin{aligned}
\Sh = 
(\zeta^4\lambda_1^{-1}&-\zeta^4)(\zeta^4\lambda_1^{-1}-\zeta^8)(\zeta^2\lambda_1^{-2}\lambda_2^{-1}-\zeta^8)(\zeta^2\lambda_1^{-2}\lambda_2^{-1}-\zeta^4)(\lambda_1^{-3}\lambda_2^{-2}+1) \\
&\times (\zeta^{10}\lambda_1^{-1}\lambda_2^{-1}-\zeta^9)(\zeta^{10}\lambda_1^{-1}\lambda_2^{-1}+1)(\zeta^{10}\lambda_1^{-1}\lambda_2^{-1}-\zeta^3)(\lambda_2^{-1}-1).
\end{aligned}
\end{align}
Then $\Sh =0$ if and only if one of the factors in \eqref{eq:shapovalov} vanishes. Let
\begin{align}\label{eq:substes-roots}
\St_1 &= \{1, \zeta^8\}, & \St_2 &= \{-1, \zeta^{10}\}, &
\St_3 &= \{\zeta, \zeta^4, \zeta^7\}.
\end{align}
The equalities alluded above can be packed as the conditions:
\begin{align}\label{eq:shapovalov-equalities}
\lambda_1 &\overset{?}{\in} \St_1,& \lambda_1^2\lambda_2 &\overset{?}{\in} \St_2,&
\lambda_1^3\lambda_2^2 &\overset{?}{=} -1,&   \lambda_1\lambda_2 &\overset{?}{\in} \St_3,&   \lambda_2 &\overset{?}{=} 1.
\end{align}
To organize the information, we consider 47 subsets of $\hgamma$, organized in classes $\cC_j$ according to the quantity $j$ of conditions in 
\eqref{eq:shapovalov-equalities} satisfied.
 The class $\cC_0$ contains just one family:
\begin{align*}
\If_1 &= \{ \lambda \in \widehat{\Gamma} \mid \lambda_1\notin \St_1, \, \lambda_1^2\lambda_2 \notin \St_2, \, \lambda_1^3\lambda_2^2 \neq -1,\,   \lambda_1\lambda_2 \notin \St_3, \,   \lambda_2 \neq 1\}; \end{align*}
Here is the class $\cC_1$:
\begin{align*}
\If_2 & = \{ \lambda \in \widehat{\Gamma} \mid \lambda_1= 1, \, \lambda_1^2\lambda_2 \notin \St_2, \, \lambda_1^3\lambda_2^2 \neq -1,\,   \lambda_1\lambda_2 \notin \St_3, \,   \lambda_2 \neq 1 \} \\
& = \{ \lambda \in \widehat{\Gamma} \mid \lambda_1= 1, \lambda_2\notin \{1,\zeta,\zeta^4,\zeta^7,\zeta^3,\zeta^9,-1,\zeta^{10} \}\}; \\
\If_3 & = \{ \lambda \in \widehat{\Gamma} \mid \lambda_1= \zeta^8,\, \lambda_1^2\lambda_2 \notin \St_2, \, \lambda_1^3\lambda_2^2 \neq -1,\,   \lambda_1\lambda_2 \notin \St_3, \,   \lambda_2 \neq 1\} \\
& = \{ \lambda \in \widehat{\Gamma} \mid \lambda_1= \zeta^8, \, \lambda_2\notin \{\pm  1, \zeta^2, \zeta^3,\zeta^5, \zeta^8,\zeta^9,\zeta^{11}\} \};
\\
\If_4 & = \{ \lambda \in \widehat{\Gamma} \mid \lambda_1\notin \St_1, \, \lambda_1^2\lambda_2=-1,\,  \lambda_1^3\lambda_2^2 \neq -1,\,   \lambda_1\lambda_2 \notin \St_3, \,   \lambda_2 \neq 1\} \\
& = \{ \lambda \in \widehat{\Gamma} \mid \lambda_1^2\lambda_2=-1,\, \lambda_1\notin \{\pm 1,\zeta^8,\zeta^{10},\zeta^4,\zeta^2\} \};
\\
\If_5 & = \{ \lambda \in \widehat{\Gamma} \mid \lambda_1\notin \St_1, \, \lambda_1^2\lambda_2=\zeta^{10},\,
\lambda_1^3\lambda_2^2 \neq -1,\,   \lambda_1\lambda_2 \notin \St_3, \,   \lambda_2 \neq 1 \} \\
& = \{ \lambda \in \widehat{\Gamma} \mid \lambda_1^2\lambda_2=\zeta^{10}, \, \lambda_1\notin \{\pm 1,\zeta^8,\zeta^{10},\zeta^4,\zeta^2\} \}; \\
\If_6 & = \{ \lambda \in \widehat{\Gamma} \mid \lambda_1\notin \St_1, \, \lambda_1^2\lambda_2 \notin \St_2, \,  \lambda_1^3\lambda_2^2=-1, \,  \lambda_1\lambda_2 \notin \St_3, \,   \lambda_2 \neq 1 \} \\
& = \{ \lambda \in \widehat{\Gamma} \mid \lambda_1^3\lambda_2^2=-1, \, \lambda_1\notin \{\pm 1,\zeta^8,\zeta^{10},\zeta^4,\zeta^2\} \}; \\
\If_7 & = \{ \lambda \in \widehat{\Gamma} \mid \lambda_1\notin \St_1, \, \lambda_1^2\lambda_2 \notin \St_2, \, \lambda_1^3\lambda_2^2 \neq -1,\,  \lambda_1\lambda_2=\zeta, \,    \lambda_2 \neq 1 \} \\
& = \{ \lambda \in \widehat{\Gamma} \mid \lambda_1\lambda_2=\zeta, \, \lambda_1\notin \{1,\zeta^8,\zeta,\zeta^4,\zeta^9\} \}; \\
\If_8 & = \{ \lambda \in \widehat{\Gamma} \mid \lambda_1\notin \St_1, \, \lambda_1^2\lambda_2 \notin \St_2, \, \lambda_1^3\lambda_2^2 \neq -1,\,   \lambda_1\lambda_2=\zeta^4, \,   \lambda_2 \neq 1 \} \\
& = \{ \lambda \in \widehat{\Gamma} \mid \lambda_1\lambda_2=\zeta^4, \, \lambda_1\notin \{1,\zeta^8,\zeta^4,\zeta^2,-1,\zeta^{10} \}\}; \\
\If_9 & = \{ \lambda \in \widehat{\Gamma} \mid \lambda_1\notin \St_1, \, \lambda_1^2\lambda_2 \notin \St_2, \, \lambda_1^3\lambda_2^2 \neq -1,\,   \lambda_1\lambda_2=\zeta^7, \, \lambda_2 \neq 1 \} \\
& = \{ \lambda \in \widehat{\Gamma} \mid \lambda_1\lambda_2=\zeta^7, \, \lambda_1\notin \{1,\zeta^8,\zeta^7,\zeta^4,\zeta^{11} \}\}; \\
\If_{10} & = \{ \lambda \in \widehat{\Gamma} \mid \lambda_1\notin \St_1, \, \lambda_1^2\lambda_2 \notin \St_2, \, \lambda_1^3\lambda_2^2 \neq -1,\,   \lambda_1\lambda_2 \notin \St_3, \,   \lambda_2=1 \}
\\ &= \{ \lambda \in \widehat{\Gamma} \mid \lambda_1\notin \G_{12}, \, \lambda_2=1 \};
\end{align*}
All the 37 remaining subsets belong to class $\cC_2$:

\begin{align*}
\If_{11} & = \{ \lambda \in \widehat{\Gamma} \mid \lambda_1= 1, \lambda_2=\zeta \}, &
\If_{12} & = \{ \lambda \in \widehat{\Gamma} \mid \lambda_1= 1, \lambda_2=\zeta^4 \},
\\
\If_{13} & = \{ \lambda \in \widehat{\Gamma} \mid \lambda_1= 1, \lambda_2=\zeta^7 \}, &
\If_{14} & = \{ \lambda \in \widehat{\Gamma} \mid \lambda_1= 1, \lambda_2=\zeta^3 \},
\\
\If_{15} & = \{ \lambda \in \widehat{\Gamma} \mid \lambda_1= 1, \lambda_2=\zeta^9 \}, &
\If_{16} & = \{ \lambda \in \widehat{\Gamma} \mid \lambda_1= 1, \lambda_2=-1 \},
\\
\If_{17} & = \{ \lambda \in \widehat{\Gamma} \mid \lambda_1= 1, \lambda_2=\zeta^{10} \}, &
\If_{18} & = \{ \lambda \in \widehat{\Gamma} \mid \lambda_1= \zeta^8 , \lambda_2=\zeta^5 \},
\\
\If_{19} & = \{ \lambda \in \widehat{\Gamma} \mid \lambda_1= \zeta^8 , \lambda_2=\zeta^8 \}, &
\If_{20} & = \{ \lambda \in \widehat{\Gamma} \mid \lambda_1= \zeta^8 , \lambda_2=\zeta^{11} \},
\\
\If_{21} & = \{ \lambda \in \widehat{\Gamma} \mid \lambda_1= \zeta^8 , \lambda_2=\zeta^3 \}, &
\If_{22} & = \{ \lambda \in \widehat{\Gamma} \mid \lambda_1= \zeta^8 , \lambda_2=\zeta^9 \},
\\
\If_{23} & = \{ \lambda \in \widehat{\Gamma} \mid \lambda_1= \zeta^8 , \lambda_2=\zeta^2 \}, &
\If_{24} & = \{ \lambda \in \widehat{\Gamma} \mid \lambda_1= \zeta^8 , \lambda_2=-1 \},
\\
\If_{25} & = \{ \lambda \in \widehat{\Gamma} \mid \lambda_1=\zeta^{11}, \lambda_2=\zeta^8 \}, &
\If_{26} & = \{ \lambda \in \widehat{\Gamma} \mid \lambda_1=\zeta^5,\lambda_2=\zeta^8 \},
\\
\If_{27} & = \{ \lambda \in \widehat{\Gamma} \mid \lambda_1=\zeta^4, \lambda_2=\zeta^9 \}, &
\If_{28} & = \{ \lambda \in \widehat{\Gamma} \mid \lambda_1=\zeta^9, \lambda_2=\zeta^4 \},
\\
\If_{29} & = \{ \lambda \in \widehat{\Gamma} \mid \lambda_1=-1, \lambda_2=-1 \}, &
\If_{30} & = \{ \lambda \in \widehat{\Gamma} \mid \lambda_1=\zeta^2, \lambda_2=\zeta^2 \},
\\
\If_{31} & = \{ \lambda \in \widehat{\Gamma} \mid \lambda_1=-1, \lambda_2=\zeta^{10} \}, &
\If_{32} & = \{ \lambda \in \widehat{\Gamma} \mid \lambda_1=\zeta^{10}, \lambda_2=-1 \},
\\
\If_{33} & = \{ \lambda \in \widehat{\Gamma} \mid \lambda_1=\zeta^2, \lambda_2=-1 \}, &
\If_{34} & = \{ \lambda \in \widehat{\Gamma} \mid \lambda_1=\zeta^4, \lambda_2=\zeta^3 \},
\\
\If_{35} & = \{ \lambda \in \widehat{\Gamma} \mid \lambda_1=\zeta^3, \lambda_2=\zeta^4 \},
\\
\If_{36} & = \{ \lambda \in \widehat{\Gamma} \mid \lambda_1=\zeta, \lambda_2=1 \}, &
\If_{37} & = \{ \lambda \in \widehat{\Gamma} \mid \lambda_1=\zeta^2, \lambda_2=1 \},
\\
\If_{38} & = \{ \lambda \in \widehat{\Gamma} \mid \lambda_1=\zeta^3, \lambda_2=1 \}, &
\If_{39} & = \{ \lambda \in \widehat{\Gamma} \mid \lambda_1=\zeta^4, \lambda_2=1 \},
\\
\If_{40} & = \{ \lambda \in \widehat{\Gamma} \mid \lambda_1=\zeta^5, \lambda_2=1 \}, &
\If_{41} & = \{ \lambda \in \widehat{\Gamma} \mid \lambda_1=-1, \lambda_2=1 \},
\\
\If_{42} & = \{ \lambda \in \widehat{\Gamma} \mid \lambda_1=\zeta^7, \lambda_2=1 \}, &
\If_{43} & = \{ \lambda \in \widehat{\Gamma} \mid \lambda_1=\zeta^8, \lambda_2=1 \},
\\
\If_{44} & = \{ \lambda \in \widehat{\Gamma} \mid \lambda_1=\zeta^9, \lambda_2=1 \}, &
\If_{45} & = \{ \lambda \in \widehat{\Gamma} \mid \lambda_1=\zeta^{10}, \lambda_2=1 \},
\\
\If_{46} & = \{ \lambda \in \widehat{\Gamma} \mid \lambda_1=\zeta^{11}, \lambda_2=1 \}, &
\If_{47} & = \{ \lambda \in \widehat{\Gamma} \mid \lambda_1=1, \lambda_2=1 \}.
\end{align*}

\bigbreak
\begin{maintheorem} \label{th:main result}
The dimension and the maximal degree of $L(\lambda)$ depend on $\lambda_i$, $i\in \I_2$,
and appear in Table \ref{tab:dim-hwv}.
\end{maintheorem}

The paper is organized as follows. We collect some general information about $\cU$ and the Verma modules in Section \ref{sect:preliminaries}, where we also deal with $\If_1$.
The proof of the Main Theorem for the families in the class 1, resp. 2, is given
in Section \ref{sec:irreps1}, respectively \ref{sec:irreps2}.

\medbreak
If $M \in \cU$, then we write $N \leq M$ to express that $N$ is a submodule of $M$.

\newpage

\begin{table}[H]
	\caption{Dimensions and highest degrees of simple modules}\label{tab:dim-hwv}
	\begin{center}
		\begin{tabular}{|c|c|c|c|}
			\hline
			Family &  $\dim L(\lambda)$  & max. degree &  $L(\lambda)^{\varphi}$ \\
			\hline
			$\If_1$  & $144$ & $(12,8)$ & $\If_{1}$
			\\ \hline
			$\If_2$  & $48$ & $(10,8)$ & $\If_{2}$
			\\ \hline
			$\If_3$ & $96$ & $(11,8)$ & $\If_{3}$
			\\ \hline
			$\If_4$  & $48$ & $(8,6)$ & $\If_{4}$
			\\ \hline
			$\If_5$ & $96$ & $(10,7)$ & $\If_{5}$
			\\ \hline
			$\If_6$ &  $72$ & $(9,6)$ & $\If_{6}$
			\\ \hline
			$\If_7$ &  $36$ & $(9,5)$ & $\If_{7}$
			\\ \hline
			$\If_8$ & $72$ & $(10,6)$ & $\If_{8}$
			\\ \hline
			$\If_9$ & $108$ & $(11,7)$ & $\If_{9}$
			\\ \hline
			$\If_{10}$ & $72$ & $(12,7)$ & $\If_{10}$
			\\ \hline
			$\If_{11}$  & $11$ & $(5,4)$ & $\If_{12}$
			\\ \hline
			$\If_{12}$  & $11$ & $(5,4)$ & $\If_{11}$
			\\ \hline
			$\If_{13}$  & $23$ & $(7,5)$ & $\If_{44}$
			\\ \hline
			$\If_{14}$  & $25$ & $(7,5)$ & $\If_{28}$
			\\ \hline
			$\If_{15}$  & $37$ & $(9,6)$ & $\If_{41}$
			\\ \hline
			$\If_{16}$  & $37$ & $(8,6)$ & $\If_{30}$
			\\ \hline
			$\If_{17}$  & $47$ & $(10,7)$ & $\If_{46}$
			\\ \hline
			$\If_{18}$ & $11$ & $(5,3)$ & $\If_{38}$
			\\ \hline
			$\If_{19}$  & $35$ & $(8,5)$ & $\If_{40}$
			\\ \hline
			$\If_{20}$  & $71$ & $(11,7)$ & $\If_{42}$
			\\ \hline
			$\If_{21}$  & $61$ & $(9,6)$ & $\If_{32}$
			\\ \hline
			$\If_{22}$  & $49$ & $(9,6)$ & $\If_{45}$
			\\ \hline
			$\If_{23}$ & $47$ & $(8,6)$ & $\If_{29}$
			\\ \hline
			$\If_{24}$ & $85$ & $(10,7)$ & $\If_{35}$
			\\
			\hline
		\end{tabular}
\quad  \begin{tabular}{|c|c|c|c|}
	\hline
	Family &  $\dim L(\lambda)$  & max. degree &  $L(\lambda)^{\varphi}$ \\
	\hline
	$\If_{25}$ &  $37$ & $(8,5)$ & $\If_{37}$
	\\ \hline
	$\If_{26}$ &  $25$ & $(8,5)$ & $\If_{43}$
	\\ \hline
	$\If_{27}$ &  $35$ & $(9,5)$ & $\If_{36}$
	\\ \hline
	$\If_{28}$ &  $25$ & $(7,5)$ & $\If_{14}$
	\\ \hline
	$\If_{29}$ &  $47$ & $(8,6)$ & $\If_{23}$
	\\ \hline
	$\If_{30}$ &  $37$ & $(8,6)$ & $\If_{16}$
	\\ \hline
	$\If_{31}$ &  $61$ & $(10,6)$ & $\If_{39}$
	\\ \hline
	$\If_{32}$ &  $61$ & $(9,6)$ & $\If_{21}$
	\\ \hline
	$\If_{33}$ &  $71$ & $(9,6)$ & $\If_{34}$
	\\ \hline
	$\If_{34}$ &   $71$ & $(9,6)$ & $\If_{33}$
	\\ \hline
	$\If_{35}$ &   $85$ & $(10,7)$ & $\If_{24}$
	\\ \hline
	$\If_{36}$ &  $35$ & $(9,5)$ & $\If_{27}$
	\\ \hline
	$\If_{37}$ &  $37$ & $(8,5)$ & $\If_{25}$
	\\ \hline
	$\If_{38}$ &  $11$ & $(5,3)$ & $\If_{18}$
	\\ \hline
	$\If_{39}$ &  $61$ & $(10,6)$ & $\If_{31}$
	\\ \hline
	$\If_{40}$ &   $35$ & $(8,5)$ & $\If_{19}$
	\\ \hline
	$\If_{41}$ &  $37$ & $(9,6)$ & $\If_{15}$
	\\ \hline
	$\If_{42}$ &  $71$ & $(11,7)$ & $\If_{20}$
	\\ \hline
	$\If_{43}$ &   $25$ & $(8,5)$ & $\If_{26}$
	\\ \hline
	$\If_{44}$ &   $23$ & $(7,5)$ & $\If_{13}$
	\\ \hline
	$\If_{45}$ &  $49$ & $(9,6)$ & $\If_{22}$
	\\ \hline
	$\If_{46}$ &  $47$ & $(10,7)$ & $\If_{17}$
	\\ \hline
	$\If_{47}$ &  $1$ & $(0,0)$ & $\If_{47}$
	\\
	\hline
\end{tabular}
	\end{center}
\end{table}

\section{Preliminaries}\label{sect:preliminaries}
\subsection{The algebra $\cU$}

The Nichols algebra $\cB(V)$ has a PBW-basis given by
\begin{align}\label{def:PBW1}
\big\{E_2^{a_2}E_{12}^{a_{12}}&E_{11212}^{a_{11212}}E_{112}^{a_{112}}E_1^{a_1} |&  a_2, a_{11212}&\in \I_{0,1};
& a_{12} &\in \I_{0,3}; &  a_{112},a_1&\in \I_{0,2}\big\}.
\end{align}
See \cite{A-unidentified}.
We obtain a new PBW-basis by reordering the PBW-generators:
\begin{align}\label{def:PBW2}
\big\{E_1^{a_1}E_{112}^{a_{112}}&E_{11212}^{a_{11212}}E_{12}^{a_{12}}E_2^{a_2} |&   a_2, a_{11212}&\in \I_{0,1};
& a_{12} &\in \I_{0,3}; & a_{112},a_1&\in \I_{0,2}\big\}.
\end{align}
Thus the set of positive roots of $\cB(V)$ (the degrees of the generators of the PBW-basis) is
\begin{align*}
\Delta_+^{V}= \left\{ \alpha_1,2\alpha_1+\alpha_2,3\alpha_1+2\alpha_2,\alpha_1+\alpha_2,\alpha_2 \right\}.
\end{align*}

By \cite[Theorem 4.9]{A-convex}, we have
\begin{align}
\label{rem:otras PRV =0}
 E_{112}^3&=E_{11212}^2=E_{12}^4=0.
\end{align}
From the defining relations \eqref{th:relations double nichols-E}, we can  deduce that
the following are valid in $\toba(V)$:
\begin{align*}
E_1E_{112}&=q_{12}\zeta^8E_{112}E_1, \\ E_{112}E_2&=-q_{12}^2E_2E_{112}+q_{12}\zeta^8E_{12}^2,
\\
E_1E_{11212}&=q_{12}^2E_{11212}E_1+q_{12}\zeta^{7}(1+\zeta)E_{112}^2
\\
E_1E_{12}^2&=E_{11212}+q_{12}\zeta(1+\zeta^3)E_{12}E_{112}+ q_{12}^2\zeta^8E_{12}^2E_1
\\  E_1E_{12}^3&= q_{12}\zeta^{10}E_{12}E_{11212}+q_{12}^2\zeta^5E_{12}^2E_{112}+ q_{12}^3E_{12}^3E_1,\\
E_1^2E_{2}&=E_{112}+q_{12}^2\zeta^2E_{12}E_1+q_{12}^2E_2E_1^2,\\ E_1^2E_{12}&=-q_{12}^2E_{112}E_1+q_{12}^2\zeta^8E_{12}E_1^2, \\
E_{112}E_{12}^2&=-q_{12}\zeta^4(1+\zeta^3)E_{12}E_{11212}+ q_{12}^2\zeta^2E_{12}^2E_{112}
\\ E_{112}E_{12}^3 &=q_{12}^2\zeta^{11}E_{12}^2E_{11212}+  q_{12}^3\zeta^3E_{12}^3E_{112},\\
E_{11212}E_{12} &=q_{12}\zeta^{10}E_{12}E_{11212},
\\ E_{112}E_{11212}&=q_{12}\zeta^{9}E_{11212}E_{112},\\
E_{11212}E_2&=q_{12}^3 E_2E_{11212}+q_{12}^2\zeta^2(1+\zeta)E_{12}^3,
\\ E_{12}E_2 &=-q_{12}E_2E_{12}.
\end{align*}

The following equalities hold by direct computation from  \eqref{th:relations double nichols-E} and the previous ones:
\begin{align*}
F_1E_{12}&=E_{12}F_1+q_{12}(\zeta-1)E_2\gpl_1^{-1},  \\
F_1E_{112}&=E_{112}F_1+q_{12}\zeta^8(1+\zeta^3)E_{12}\gpl_1^{-1},  \\
F_1E_{11212}&=E_{11212}F_1+q_{12}^2(\zeta^5-1)E_{12}^2\gpl_1^{-1},  \\
F_1E_{112}^2&=E_{112}^2F_1-q_{12}(1+\zeta^3)(E_{11212}\gpl_1^{-1}+ \zeta^4E_{112}E_{12}\gpl_1^{-1}),\\
F_1E_{12}^2&=E_{12}^2F_1+q_{12}^2(3)_{\zeta^5}E_2E_{12}\gpl_1^{-1},
\\
F_1E_{12}^3&=E_{12}^2F_1+q_{12}^3\zeta^3(\zeta-1)E_2E_{12}^2\gpl_1^{-1},\\
F_2E_{12}&=E_{12}F_2+(\zeta^{11}-1)E_1 g_2, \\
F_2E_{112}&=E_{112}F_2-(3)_{\zeta^7}E_1^2g_2, \\
F_2E_{11212}&=E_{11212}F_2-E_{112}E_1g_2, \\
F_2E_{12}^2&=E_{12}^2F_2+q_{21}(1+\zeta^{5})E_{112}g_2-(3)_{\zeta^{7}}E_{12}E_{1}g_2,\\
F_2E_{112}^2&=E_{112}^2F_2+(3)_{\zeta^7}\zeta^4E_{112}E_1^2g_2, \\
F_2E_{12}^3&=E_{12}^3F_2+\zeta^8(1-\zeta)(E_{12}^2E_1g_2- q_{21}\zeta^{3}E_{12}E_{112}g_2+q_{21}^2\zeta^{3}E_{11212}g_2),
\\
F_{11212}E_{11212}&=E_{11212}F_{11212}+\gpl_1^{-3}\gpl_2^{-2}-g_{11212},
 \\
F_{12} E_2 &= E_2 F_{12} +(1-\zeta^{11})F_1\gpl^{-1}_2,  \\
F_{12}E_{12}&=E_{12}F_{12}+\gpl_1^{-1}\gpl_2^{-1}-g_1g_2, \\
F_{12} E_{112}&= E_{112} F_{12}+\zeta^3(3)_{\zeta^7}E_1 g_1g_2,\\
F_{12} E_{112}^2&= E_{112}^2 F_{12}+\zeta^{11}(3)_{\zeta^7}E_{112}E_1 g_1g_2,
\\
 F_{12} E_1&= E_1 F_{12}+q_{21}(1-\zeta)F_2 g_1, \\
F_{12} E_{11212}&= E_{11212} F_{12}+\zeta^{11}E_{112} g_1g_2,\\
F_{112}E_{112}&=E_{112}F_{112}+\gpl_1^{-2}\gpl_2^{-1}-g_1^2g_2,\\
 F_{112}E_2&= E_2 F_{112} + (\zeta-1)F_1^2\gpl^{-1}_2.
\end{align*}

\subsection{Verma modules}\label{subsec:verma-gral}\label{subsec:reps}

We shall use the notation for $q$-factorial numbers: for each $q\in\ku^\times$,
\begin{align*}
(n)_q&=1+q+\ldots+q^{n-1}, & (n)_q!&=(1)_q(2)_q\cdots(n)_q,& n&\in\N.
\end{align*}

We shall investigate the lattice of submodules of a Verma module.
We record the following standard fact for future use.

\begin{remark}\label{rem:F actua por 0}
	Let $v\in M(\lambda)_\alpha$ be such that $F_i\cdot v=0$ for $i\in \I_2$. By the triangular decomposition of $\Uc$, $\Uc\cdot v=\Uc^+\cdot v$.
	In particular, if  $\alpha\neq0$, then $\Uc\cdot v \cap \ku v_\lambda= 0$.
\end{remark}

We  consider two families in $M(\lambda)$,
corresponding to  PBW-bases \eqref{def:PBW1} and \eqref{def:PBW2}. We set
\begin{align*}
\widetilde{m}_{a,b,c,d,e}&:=E_2^aE_{12}^bE_{11212}^cE_{112}^dE_1^e\cdot v_\lambda, & \widetilde{n}_{a,b,c,d,e}&:=E_1^e E_{112}^d E_{11212}^c E_{12}^b E_2^a\cdot v_\lambda
\end{align*}
for $a,b,c,d,e \in\Z$. Clearly, $v_\lambda=\widetilde{m}_{0,0,0,0,0}=\widetilde{n}_{0,0,0,0,0}$ and

\begin{align*}
\widetilde{m}_{a,b,c,d,e}\neq0   \iff
a,c \in \I_{0,1}, b \in\I_{0,3}, d,e \in \I_{0,2} \iff \widetilde{n}_{a,b,c,d,e}\neq0.
\end{align*}

\smallbreak
We denote by $\langle S \rangle$ the subspace generated by a subset $S$ of a vector space. Let
\begin{align*}
W_1(\lambda)&=\langle\widetilde{m}_{a,b,c,d,e} \mid a,c \in \I_{0,1}, b \in\I_{0,3}, d \in \I_{0,2}, e \in \I_{1,2}\rangle,\\
W_2(\lambda)&=\langle\widetilde{m}_{a,b,c,d,2} \mid   a,c \in \I_{0,1}, b \in\I_{0,3}, d \in \I_{0,2}  \rangle,\\
W(\lambda)&=\langle\widetilde{n}_{1,b,c,d,e}\mid c \in \I_{0,1}, b \in\I_{0,3}, d,e \in \I_{0,2}\rangle.
\end{align*}

By a direct computation, we can prove:

\begin{lemma}\label{coro:Fi actuando en Wj}
\begin{enumerate}[leftmargin=*,label=\rm{(\alph*)}]
\item\label{(a)}  $F_2 \cdot W_i(\lambda) \subseteq  W_i(\lambda)$, $i\in \I_2$,
\item $  F_1 \cdot \widetilde{m}_{a,b,c,d,i} \in  \lambda(\gpl_1^{-1})(i)_{\zeta^4}(\zeta^{(i-1)8}-\lambda_1)\widetilde{m}_{a,b,c,d,i-1} + W_i(\lambda)$,
$i\in \I_2$,
\item\label{(b)}  $F_1 \cdot W(\lambda) \subseteq  W(\lambda)$,
\item
$F_2 \cdot \widetilde{n}_{1,b,c,d,e} \in  \lambda(\gpl_2^{-1})(1-\lambda_2)\widetilde{n}_{0,b,c,d,e} + W(\lambda)$.
\end{enumerate}

\medbreak
In consequence,
\begin{itemize}[leftmargin=*]\renewcommand{\labelitemi}{$\circ$}
 \smallbreak \item $W_1(\lambda)$ is a $\Uc$-submodule if and only if $\lambda_1=1$;

 \smallbreak \item $W_2(\lambda)$ is a $\Uc$-submodule if and only if $\lambda_1=\zeta^8$;

\smallbreak \item $W(\lambda)$ is a $\Uc$-submodule if and only if $\lambda_2=1$. \qed
\end{itemize}
\end{lemma}

We denote by $m_{a,b,c,d,e}$,  $n_{a,b,c,d,e}$ the classes of $\widetilde{m}_{a,b,c,d,e}$,  $\widetilde{n}_{a,b,c,d,e}$ in $L(\lambda)$.
We order lexicographically the set of all  $m_{a,b,c,d,e}$: 
\begin{align}\label{eq:lexico}
m_{a,b,c,d,e} < m_{a',b',c',d',e'} \iff a<a', \text{ or } a=a', b<b', \text{ or \dots}.
\end{align}

\subsection{Simple modules}\label{subsec:simple-particular}

Let $\varphi: \Uc \rightarrow \Uc$ be the algebra automorphism such that
\begin{align*}
\varphi(K_i) &= K_i^{-1},& \varphi(L_i)&=L_i^{-1},& \varphi(E_i) &= F_iL_i^{-1},& \varphi(F_i) &= K_i^{-1}E_i,
\end{align*}
$i\in \I_2$, cf. \cite[Proposition 4.9]{H-lusztig iso}; this resembles the Chevalley involution.
If $M$ is a $\Uc$-module, then we denote by $M^{\varphi}$ the $\cU$-module with $M^{\varphi}=M$ as vector space and action given by $a\rhd v= \varphi(a) \cdot v$, $v\in V$, $a\in \cU$.
If $v \in M$ has weight $\lambda$ (with respect the action of $\Gamma$), then $v \in M^{\varphi}$ has weight $\lambda^{-1}$.
The functor $M\mapsto M^{\varphi}$ preserves simple objects and sends lowest weight modules to highest weight modules, and vice versa. The following result is standard.

\begin{lemma}\label{lem: irreduble via phi}
The subspace $X(\lambda):=\{x\in L(\lambda): E_ix=0\mbox{ for all }i\}$ of $L(\lambda)$ is one-dimensional  and there exists
$\mu\in\widehat\Gamma$ such that $X(\lambda) \overset{(1)}{=} L(\lambda)_\mu$, $L(\lambda)^\varphi \overset{(2)}{\simeq} L(\mu^{-1})$.
\end{lemma}
\pf
$X(\lambda) \neq 0$ because there exists $\beta\in\N_0^2$ maximal such that $L(\lambda)_\beta\neq 0$.
Since $X(\lambda)$ is $\Gamma$-stable, there exists a weight vector
$0 \neq x\in X(\lambda)$ with weight $\mu\in\widehat\Gamma$. Thus $\Uc^- x= \Uc x= L(\lambda)$ and (1) follows. Also 
$L(\lambda)^\varphi = (\Uc^- x)^\varphi\twoheadrightarrow  L(\mu^{-1})$ implying (2).
\epf

\begin{lemma}\label{lem: morphism between submodules}
Let $M \in \umod$ a highest weight module of highest weight $\mu$ and  $0 \neq v\in M^{\mu}$. 
If  $m_{a,b,c,d,e}\neq 0$ in $L(\mu^{-1})$ then
$z := F_2^aF_{12}^bF_{11212}^cF_{112}^dF_1^e v \neq 0$.
\end{lemma}

There is an analogue statement for $n_{a,b,c,d,e}$.

\pf
Indeed $M^{\varphi}$ is  lowest weight of lowest weight $\mu^{-1}$, hence  
$M^{\varphi}\twoheadrightarrow L (\mu^{-1})$; up to a non-zero scalar, $z \mapsto m_{a,b,c,d,e} \neq 0$,
hence $z \neq 0$.
\epf

\subsection{A relative of $\mathfrak{u}_q(\mathfrak{sl}_2)$}\label{subsec:dim V=1}

We consider for a moment the algebra $\cV$ constructed as $\cU$ above but starting from a braided vector space of dimension $1$, with braiding given by
$q=\gpl(g) \in \G_N'$,  $g\in \Lambda$, $\gpl\in\widehat{\Lambda}$.
The algebra $\cV$ is close to $\mathfrak{u}_q(\mathfrak{sl}_2)$ and
has a presentation by generators $h\in \Lambda$, $\tau \in\widehat\Lambda$, $E$, $F$ with relations
\begin{align*}
E^N=F^N&=0, & h E&= \gpl(h) E h, & \tau E&= \tau(g) E\tau,\\
EF-FE &= g-\gpl^{-1}, & h F&= \gpl^{-1}(h)F h, & \tau F&= \tau(g^{-1}) F\tau,
\end{align*}
and $h\tau=\tau h$ for $h\in\Lambda$, $\tau\in\widehat\Lambda$, and the relations defining $\Lambda$, $\widehat\Lambda$. Thus
\begin{align}\label{eq:commute Ej F}
E^jF-FE^j&= (j)_q E^{j-1}(g-q^{1-j}\gpl^{-1}), & j &\in \N.
\end{align}

Let $\lambda \in \widehat{\Gamma}$. Let $L(\lambda)$ be lowest weight $\cV$-module of lowest weight $\lambda$ defined in the same usual way.
The same argument as for $\mathfrak{u}_q(\mathfrak{sl}_2)$ gives the following.

\begin{lemma}\label{coro:modulos uqsl2}
\begin{enumerate}[leftmargin=*,label=\rm{(\alph*)}]
  \item If there exists $j\in\I_{N-1}$ such that $\lambda(g\gpl)=q^{1-j}$, then $\dim L(\lambda)=j$.

  \smallbreak
  \item If $\lambda(g\gpl)\notin\{q^h|h \in\I_{0,N-2}\}$, then $\dim L(\lambda)=N$.

  \smallbreak
  \item $L(\lambda)$ has a basis $v_0,\dots,v_{\dim L(\lambda)-1}$ such that for all $i$,
\begin{align}
Ev_i&=v_{i+1}, & Fv_i&=(i)_q(q^{1-i}\lambda(\gpl_1^{-1})-\lambda(g_1))v_{i-1}, & h\tau v_i&=\lambda(h\tau) \gpl^i(h)\tau(g^i) v_i. 
\end{align}

  \smallbreak
  \item
Let $M$ be a lowest weight $\cV$-module with lowest weight $\lambda\in\widehat{\Gamma}$. If $0\neq v\in M^\lambda$, 
then $v, E v, \dots, E^{n-1}v$ are linearly independent, where
\begin{enumerate}[label=\rm{(\arabic*)}]

\smallbreak  \item either $n=j$ if $\lambda(g\gpl)=q^{1-j}$ for some (unique) $j\in\I_{N-1}$,

\smallbreak  \item or else $n=N-1$ if $\lambda(g\gpl)\notin\{q^h|h \in\I_{0,N-2}\}$.
\end{enumerate}

\smallbreak 
Moreover $F^i E^i v=a_i v$ for some $a_i\in\ku^\times$ when $i\in\I_{0, n-1}$. \qed
\end{enumerate}
\end{lemma}

\subsection{The class $\cC_0$}

The first family is easy to deal with.

\begin{lemma}\label{lem:modulo irr caso 1}
If $\lambda \in \If_1$, then $M(\lambda)$ is simple.
\end{lemma}
\pf
By \cite[5.16]{HY} that says:  if 
$\Sh \neq 0$, 
then $M(\lambda)$ is simple.
\epf

\section{Simple $\Uc$-modules in class $\cC_1$}\label{sec:irreps1}
Here we deal with the class of families satisfying exactly one of the conditions in \eqref{eq:shapovalov-equalities}.
Recall that $\Gamma = \Lambda \times \VLambda$; we introduce $\chi_i \in \VGamma$ by
\begin{align*}
\chi_i (g, \gpl) &= \gpl_i(g) \gpl(g_i),& i &\in \I_2.
\end{align*}
For simplicity, we introduce the following notation:
\begin{align*}
g_{12} &= g_1g_2, & g_{112} &= g_1^2g_2, & g_{11212} &= g_1^3g_2^2, \\ 
\gpl_{12} &= \gpl_1\gpl_2, & \gpl_{112} &= \gpl_1^2\gpl_2, & \gpl_{11212} &= \gpl_1^3\gpl_2^2.
\end{align*}

We outline the method to compute $L(\lambda)$, $\lambda\in\If_j$, $j\in\I_{2,10}$.

\begin{enumerate}[leftmargin=*, label=\rm{(\alph*)}]
\item As (exactly) one of the factors of the Shapovalov determinant $\Sh$ vanishes, there exists $\beta \neq 0$ and 
$w\in M(\lambda)_\beta - 0$, such that $F_i w=0$, $i\in\I_2$, 
see Remarks \ref{rem:rel E112}, \ref{rem:rel E112 cuadrado}, \ref{rem:l1 a la 3 l2 a la 2}, \ref{rem:rel E12}, \ref{rem:rel E12 cuadrado}, \ref{rem:rel E12 cubo},
or Lemma \ref{coro:Fi actuando en Wj}.
Thus $\Uc w$ is a proper submodule.

\smallbreak
\item Assume we are dealing with $\If_j$, $j\in \I_{2,6}$. Write $w=\sum \mathtt{p}_{a,b,c,d,e} \ \widetilde{m}_{a,b,c,d,e}$. 
Then there exists $a,b,c,d,e$ such that $\mathtt{p}_{a,b,c,d,e} \neq 0$ and exactly four of the integers $a, \dots, e$ are zero. 
The same holds for $j\in\I_{7,10}$ exchanging  $\widetilde{m}_{a,b,c,d,e}$ by $\widetilde{n}_{a,b,c,d,e}$.
From here we describe a basis $\mathrm{B}_j$ of the quotient $L'(\lambda)$  of $M(\lambda)$ by $\Uc w$, $j\in\I_{2,10}$.

\smallbreak
\item Let $v$ be the element of maximal degree of $L'(\lambda)$. A short computation shows that $v$ belongs to every submodule of $L'(\lambda)$.
Because of the inequalities defining $\If_j$, there exists $F\in\Uc$ such that $Fv=v_\lambda$. Hence $L'(\lambda)$ is simple.
\end{enumerate}

We work out the details for $\If_2$, with shorter expositions for the other families in $\cC_1$.

\subsection{The family $\If_2$} Recall that
\begin{align*}
\If_2 &  = \{ \lambda \in \widehat{\Gamma} \mid \lambda_1= 1, \lambda_2\notin \{1,\zeta,\zeta^4,\zeta^7,\zeta^3,\zeta^9,-1,\zeta^{10} \}\}.
\end{align*}

\begin{lemma}\label{lem:modulo irr caso 2}
If $\lambda \in \If_2$, then $\dim L(\lambda)=48$. A basis of $L(\lambda)$ is given by
\begin{align*}
\mathrm{B}_2=\{m_{a,b,c,d,0}: a,c \in \I_{0,1},\, b \in\I_{0,3}, \, d \in \I_{0,2}\}.
\end{align*}
\end{lemma}
\pf
Let $w = \widetilde{m}_{0,0,0,0,1}$; then $F_i w =0$, $i\in \I_2$, 
hence $\Uc^+w = W_1(\lambda) \le M(\lambda)$ is  proper by Lemma \ref{coro:Fi actuando en Wj}. Let  
$L'(\lambda)=M(\lambda)/\Uc^+ w$. Let $\widehat{m}_{a,b,c,d,0}$ be the class of $\widetilde{m}_{a,b,c,d,0}$ in $L'(\lambda)$.
Then $$\widehat{\mathrm{B}}_2 = \{\widehat{m}_{a,b,c,d,0}: a,c \in \I_{0,1},\, b \in\I_{0,3}, \, d \in \I_{0,2}\}$$
is a basis of $L'(\lambda)$, ordered by \eqref{eq:lexico}. Thus, it is enough to show that $L'(\lambda)$ is simple. 
Let $0\neq W \leq L'(\lambda)$ and pick $u\in W-0$. Fix  $\widehat{m}_{a,b,c,d,0} \in \widehat{\mathrm{B}}_2$ minimal   among
those whose coefficient in $u$ is non-zero.  
Then
\begin{align*}
E_{112}^{2-d}E_{11212}^{1-c}E_{12}^{3-b}E_2^{1-a}  u &\in \kut  \widehat{m}_{1,3,1,2,0}  \implies \widehat{m}_{1,3,1,2,0}\in W.
\end{align*}

By abuse of notation, we denote by $v_{\lambda}$ its class  in $L'(\lambda)$.
We claim that 
\begin{align}\label{eq:If2-claim}
F_2 F_{12}^3F_{11212}F_{112}^2  \widehat{m}_{1,3,1,2,0}\in \kut v_\lambda;
\end{align}
 this implies that $v_\lambda\in W$, 
so $L'(\lambda)$  is simple.

To prove \eqref{eq:If2-claim}, we first consider the subalgebra $\cV_1 = \ku\langle g,\gpl,E_{112}, F_{112} \rangle$ of $\cU$;
clearly $\cV_1 \simeq \cV$ from \S \ref{subsec:dim V=1}. Then 
\begin{align*}
F_{112}  \widehat{m}_{1,3,1,0,0}&=0, & g_{112}\gpl_{112} \widehat{m}_{1,3,1,0,0}&= -\lambda_2 \widehat{m}_{1,3,1,0,0}, &
E_{112}^2  \widehat{m}_{1,3,1,0,0} &= \gpl_{112}^2(g_{12}^{-6}) \widehat{m}_{1,3,1,2,0}.
\end{align*}
By Lemma \ref{coro:modulos uqsl2}, we conclude that
\begin{align*}
F_{112}^2  \widehat{m}_{1,3,1,2,0} &\in \kut \widehat{m}_{1,3,1,0,0} \implies \widehat{m}_{1,3,1,0,0} \in W.
\end{align*}
We next consider $\cV_2 = \ku\langle g,\gpl,E_{11212}, F_{11212} \rangle \hookrightarrow  \cU$;
again, $\cV_2 \simeq \cV$. Then 
\begin{align*}
F_{11212}  \widehat{m}_{1,3,0,0,0}&=0, & g_{11212}\gpl_{11212}  \widehat{m}_{1,3,0,0,0}&= -\lambda_2^2 \widehat{m}_{1,3,0,0,0}, \\
E_{11212}  \widehat{m}_{1,3,0,0,0} &= \gpl_{11212}(g_1^{-3}g_2^{-4}) \widehat{m}_{1,3,1,0,0},& &\\
\overset{\text{Lemma \ref{coro:modulos uqsl2}}}{\implies} F_{11212}  \widehat{m}_{1,3,1,0,0} &\in \kut \widehat{m}_{1,3,0,0,0}& \implies \widehat{m}_{1,3,0,0,0} &\in W.
\end{align*}
Once again, we consider $\cV_3 = \ku\langle g,\gpl,E_{12}, F_{12} \rangle \hookrightarrow  \cU$;
thus $\cV_3 \simeq \cV$ from \S \ref{subsec:dim V=1}. Then
\begin{align*}
&\begin{aligned}
F_{12}  \widehat{m}_{1,0,0,0,0}&=0, & g_{12}\gpl_{12}  \widehat{m}_{1,0,0,0,0}&= \lambda_2\zeta^{11} \widehat{m}_{1,0,0,0,0}, &
E_{12}^3  \widehat{m}_{1,0,0,0,0} &= \gpl_{12}^3(g_2^{-1}) \widehat{m}_{1,3,0,0,0}
\end{aligned}
\\
&\overset{\text{Lemma \ref{coro:modulos uqsl2}}}{\implies} F_{12}^3 \widehat{m}_{1,3,0,0,0} 
\in \kut \widehat{m}_{1,0,0,0,0} \implies \widehat{m}_{1,0,0,0,0} \in W.
\end{align*}
Now $F_2\widehat{m}_{1,0,0,0,0} = \lambda(\sigma_2)^{-1}(\lambda_2 - 1) v_{\lambda} \neq 0$, and 
\eqref{eq:If2-claim} follows. 
\epf

\begin{coro}\label{cor:modulo irr caso 2}
If $\lambda \in \If_2$, then $N(\lambda) \simeq L(\chi_1\lambda)$ and $\chi_1\lambda \in \If_3$. 
\end{coro}

\pf By the proof of the Lemma, $N(\lambda)$ is of lowest weight  $\chi_1\lambda$ and $\dim N(\lambda) = 96$.
It is easy to see that $\chi_1\lambda \in \If_3$; hence $\dim L(\chi_1\lambda)=96$ by Lemma \ref{lem:modulo irr caso 3} and the claim follows.
\epf

\subsection{The family $\If_3$} Recall that
\begin{align*}
\If_3 &  =  \{ \lambda \in \widehat{\Gamma} \mid \lambda_1= \zeta^8, \, \lambda_2\notin \{\pm  1, \zeta^2, \zeta^3,\zeta^5, \zeta^8,\zeta^9,\zeta^{11}\} \}.
\end{align*}

\begin{lemma}\label{lem:modulo irr caso 3}
If $\lambda \in \If_3$, then $\dim L(\lambda)=96$. A basis of $L(\lambda)$ is given by
\begin{align*}
\mathrm{B}_3=\{m_{a,b,c,d,e}|a,c \in \I_{0,1},\, b \in\I_{0,3}, \, d \in \I_{0,2}, e\in \I_{0,1}\}.
\end{align*}
\end{lemma}

\pf Let $w = \widetilde{m}_{0,0,0,0,2}$ and   
$L'(\lambda)=M(\lambda)/\Uc^+ w$. We identify $\mathrm{B}_3$ with a basis of $L'(\lambda)$. 
Now $F_2 F_{12}^3F_{11212}F_{112}^2F_1  m_{1,3,1,2,1}\in \ku^{\times}v_\lambda$, hence $L'(\lambda)$ is simple.
\epf

Exactly as for Corollary \ref{cor:modulo irr caso 2}, we conclude:

\begin{coro}\label{cor:modulo irr caso 3}
	If $\lambda \in \If_3$, then $N(\lambda) \simeq L(\chi_1^2\lambda)$ and $\chi_1^2\lambda \in \If_2$. \qed
\end{coro}

\subsection{The family $\If_4$} Recall that
\begin{align*}
\If_4 &  = \{ \lambda \in \widehat{\Gamma} \mid \lambda_1^2\lambda_2=-1,\, \lambda_1\notin \{\pm 1,\zeta^8,\zeta^{10},\zeta^4,\zeta^2\} \}.
\end{align*}

We start by a Remark that will be useful elsewhere. 

\begin{remark}\label{rem:rel E112} Let $\lambda \in \widehat{\Gamma}$.
If $\lambda_1^2\lambda_2=-1$, then $w=F_1^2E_{112}E_1^2 v_\lambda \in M(\lambda)$ satisfies  
\begin{align}\label{eq:Fi=0}
F_1w=F_2w=0.
\end{align}
\end{remark}	

\pf	By a direct computation,
	\begin{align*}
	F_{112}E_{112}E_1^2 v_\lambda &= \lambda(\gpl_1^{-2}\gpl_2^{-1})q_{21}^2\zeta^4(\lambda_1^2\lambda_2+1)E_1^2v_\lambda.
	\end{align*}
	As $M(\lambda)_{4\alpha_1}=M(\lambda)_{3\alpha_1}=0$, we have that $F_2 E_{112}E_1^2v_\lambda=F_1 E_{112}E_1^2v_\lambda=0$, so
	\begin{align*}
	0&=F_{112}E_{112}E_1^2 v_\lambda = \zeta^8q_{12}^2 F_2F_1^2 E_{112}E_1^2 v_\lambda.
	\end{align*}
	This shows that $F_2w=0$; on the other hand, $F_1w= F_1^3(E_{112}E_1^2 v_\lambda) =0$, since $F_1^3=0$.
\epf

\begin{lemma}\label{lem:modulo irr caso 4}
If $\lambda \in \If_4$, then $\dim L(\lambda)=48$. A basis of $L(\lambda)$ is given by
\begin{align*}
\mathrm{B}_4=\{m_{a,b,c,0,e}: a,c \in \I_{0,1}, b \in\I_{0,3}, e \in \I_{0,2}\}.
\end{align*}
\end{lemma}
\pf
Let $w=F_1^2E_{112}E_1^2 v_\lambda$. By Remark \ref{rem:rel E112}, $\Uc w$ is a proper submodule. 
We identify $\mathrm{B}_4$ with a basis of $L'(\lambda):= M(\lambda)/\Uc w$. 
We check that there exists $F\in\Uc$ such that $F m_{1,3,1,0,2}=v_\lambda$. Then $L'(\lambda)$ is simple.
\epf

Exactly as for Corollary \ref{cor:modulo irr caso 2}, we conclude:

\begin{coro}\label{cor:modulo irr caso 4}
	If $\lambda \in \If_4$, then $N(\lambda) \simeq L(\chi_1^2\chi_2\lambda)$ and $\chi_1^2\chi_2\lambda \in \If_5$. 
\qed \end{coro}

\subsection{The family $\If_5$} Recall that
\begin{align*}
\If_5 &  = \{ \lambda \in \widehat{\Gamma} \mid \lambda_1^2\lambda_2=\zeta^{10}, \, \lambda_1\notin \{\pm 1,\zeta^8,\zeta^{10},\zeta^4,\zeta^2\} \}.
\end{align*}

Here is another Remark that will be useful later, proved as Remark \ref{rem:rel E112}. 

\begin{remark}\label{rem:rel E112 cuadrado} Let $\lambda \in \widehat{\Gamma}$.
If $\lambda_1^2\lambda_2=\zeta^{10}$, then $w=F_1^2E_{112}^2E_1^2 v_\lambda \in M(\lambda)$ satisfies  \eqref{eq:Fi=0}.
\end{remark}

\begin{lemma}\label{lem:modulo irr caso 5}
If $\lambda \in \If_5$, then $\dim L(\lambda)=96$.
A basis of $L(\lambda)$ is given by
\begin{align*}
\mathrm{B}_5=\{m_{a,b,c,d,e}|a,c, d \in \I_{0,1}, b \in\I_{0,3},  e \in \I_{0,2} \}.
\end{align*}
\end{lemma}
\pf
Let $w=F_1^2E_{112}^2E_1^2 v_\lambda$. By Remark \ref{rem:rel E112 cuadrado}, $\Uc w$ is a proper submodule. 
We identify $\mathrm{B}_5$ with a basis of $L'(\lambda):= M(\lambda)/\Uc w$. 
We check that there exists $F\in\Uc$ such that $F m_{1,3,1,1,2}=v_\lambda$. Then $L'(\lambda)$ is simple.
\epf

Exactly as for Corollary \ref{cor:modulo irr caso 2}, we conclude:

\begin{coro}\label{cor:modulo irr caso 5}
	If $\lambda \in \If_5$, then $N(\lambda) \simeq L(\chi_1^4\chi_2^2\lambda)$ and $\chi_1^4\chi_2^2\lambda \in \If_4$. 
\qed \end{coro}

\subsection{The family $\If_6$} Recall that
\begin{align*}
\If_6 &  = \{ \lambda \in \widehat{\Gamma} \mid \lambda_1^3\lambda_2^2=-1, \, \lambda_1\notin \{\pm 1,\zeta^8,\zeta^{10},\zeta^4,\zeta^2\} \}.
\end{align*}

Still another  Remark useful elsewhere, with an analogous proof as above.

\begin{remark}\label{rem:l1 a la 3 l2 a la 2} Let $\lambda \in \widehat{\Gamma}$.
	If $\lambda_1^3\lambda_2^2=-1$, then $w=F_1^2F_{112}^2E_{11212}E_{112}^2E_1^2v_\lambda$ satisfies  \eqref{eq:Fi=0}.
\end{remark}

\begin{lemma}\label{lem:modulo irr caso 6}
If $\lambda \in \If_6$, then $\dim L(\lambda)=72$. A basis of $L(\lambda)$ is given by
\begin{align*}
\mathrm{B}_6=\{m_{a,b,0,d,e} | a\in \I_{0,1}, b \in\I_{0,3}, d,e \in \I_{0,2} \}.
\end{align*}
\end{lemma}
\pf
Let $w$ be as in Remark \ref{rem:l1 a la 3 l2 a la 2}; then  $\Uc w$ is proper. 
Again $\mathrm{B}_6$ is identified with a basis of $L'(\lambda)=M(\lambda)/\Uc w$; since there is $F\in\Uc$ 
such that $Fm_{1,3,0,2,2}=v_\lambda$, $L'(\lambda)$  is simple.
\epf

Exactly as for Corollary \ref{cor:modulo irr caso 2}, we conclude:

\begin{coro}\label{cor:modulo irr caso 6}
	If $\lambda \in \If_6$, then $N(\lambda) \simeq L(\chi_1^3\chi_2^2\lambda)$ and $\chi_1^3\chi_2^2\lambda \in \If_6$. 
\qed \end{coro}

\subsection{The family $\If_7$} Recall that
\begin{align*}
\If_7 &  = \{ \lambda \in \widehat{\Gamma} \mid \lambda_1\lambda_2=\zeta, \, \lambda_1\notin \{1,\zeta^8,\zeta,\zeta^4,\zeta^9\} \}.
\end{align*}

Again we start by a useful remark.

\begin{remark}\label{rem:rel E12} Let $\lambda \in \widehat{\Gamma}$.
If $\lambda_1\lambda_2=\zeta$, then $w=F_2E_2E_{12}v_\lambda \in M(\lambda)$ satisfies \eqref{eq:Fi=0}.
\end{remark}

\begin{lemma}\label{lem:modulo irr caso 7}
If $\lambda \in \If_7$, then $\dim L(\lambda)=36$. A basis of $L(\lambda)$ is given by
\begin{align*}
\mathrm{B}_7=\{n_{a,0,c,d,e} | a,c \in \I_{0,1}, d,e \in \I_{0,2} \}.
\end{align*}
\end{lemma}
\pf
Let $w=F_2E_2E_{12}v_\lambda$.
By Remark \ref{rem:rel E12}, $\Uc w \subsetneq M(\lambda)$. Let $L'(\lambda)=M(\lambda)/\Uc w$, so $B_7$ is a basis of $L'(\lambda)$. There exists $F\in\Uc$ such that $Fn_{1,0,1,2,2}=v_\lambda$. Then $L'(\lambda)$ is simple.
\epf

Exactly as for Corollary \ref{cor:modulo irr caso 2}, we conclude:

\begin{coro}\label{cor:modulo irr caso 7}
	If $\lambda \in \If_7$, then $N(\lambda) \simeq L(\chi_1\chi_2\lambda)$ and $\chi_1\chi_2\lambda \in \If_9$. 
\qed \end{coro}

\subsection{The family $\If_8$} Recall that
\begin{align*}
\If_8 &  = \{ \lambda \in \widehat{\Gamma} \mid \lambda_1\lambda_2=\zeta^4, \, \lambda_1\notin \{1,\zeta^8,\zeta^4,\zeta^2,-1,\zeta^{10} \} \}.
\end{align*}

\begin{remark}\label{rem:rel E12 cuadrado}Let $\lambda \in \widehat{\Gamma}$.
	If $\lambda_1\lambda_2=\zeta^4$, then $w=F_2E_2E_{12}^2 v_\lambda \in M(\lambda)$ satisfies  \eqref{eq:Fi=0}.
\end{remark}	

\pf	
Analogous to Remark \ref{rem:rel E112}.
\epf

\begin{lemma}\label{lem:modulo irr caso 8}
If $\lambda \in \If_8$, then $\dim L(\lambda)=72$. A basis of $L(\lambda)$ is given by
\begin{align*}
\mathrm{B}_8=\{n_{a,b,c,d,e}| a,b,c \in \I_{0,1}, d,e \in \I_{0,2} \}.
\end{align*}
\end{lemma}
\pf
Let $w=F_2E_2E_{12}^2v_\lambda$. By Remark \ref{rem:rel E12 cuadrado}, $\Uc w \subsetneq M(\lambda)$. Now $\mathrm{B}_8$ identifies with a basis of $L'(\lambda) := M(\lambda)/\Uc w$. Since there is $F\in\Uc$ such that $Fn_{1,1,1,2,2}=v_\lambda$, $L'(\lambda)$ is simple.
\epf

Exactly as for Corollary \ref{cor:modulo irr caso 2}, we conclude:

\begin{coro}\label{cor:modulo irr caso 8}
	If $\lambda \in \If_8$, then $N(\lambda) \simeq L(\chi_1^2\chi_2^2\lambda)$ and $\chi_1^2\chi_2^2\lambda \in \If_8$. 
\qed \end{coro}

\subsection{The family $\If_9$} Recall that
\begin{align*}
\If_9 &  = \{ \lambda \in \widehat{\Gamma} \mid \lambda_1\lambda_2=\zeta^7, \, \lambda_1\notin \{1,\zeta^8,\zeta^7,\zeta^4,\zeta^{11} \}\}.
\end{align*}

\begin{remark}\label{rem:rel E12 cubo}
Let $\lambda \in \widehat{\Gamma}$. If $\lambda_1\lambda_2=\zeta^7$, then $w=F_2E_2E_{12}^3 v_\lambda \in M(\lambda)$ satisfies  \eqref{eq:Fi=0}.
\end{remark}	

\pf	
Analogous to Remark \ref{rem:rel E112}.
\epf

\begin{lemma}\label{lem:modulo irr caso 9}
If $\lambda \in \If_9$, then $\dim L(\lambda)=108$. A basis of $L(\lambda)$ is given by
\begin{align*}
\mathrm{B}_9=\{n_{a,b,c,d,e}| a,c \in \I_{0,1}, b,d,e \in \I_{0,2} \}.
\end{align*}
\end{lemma}
\pf
Let $w=F_2E_2E_{12}^3v_\lambda$. By Remark \ref{rem:rel E12 cubo}, $\Uc w \subsetneq M(\lambda)$. Let $L'(\lambda)=M(\lambda)/\Uc w$, so $\mathrm{B}_9$ is a basis of $L'(\lambda)$. Since there exists $F\in\Uc$ such that $Fn_{1,2,1,2,2}=v_\lambda$, $L'(\lambda)$ is simple.
\epf

Exactly as for Corollary \ref{cor:modulo irr caso 2}, we conclude:

\begin{coro}\label{cor:modulo irr caso 9}
If $\lambda \in \If_9$, then $N(\lambda) \simeq L(\chi_1^3 \chi_2^3 \lambda)$ and $\chi_1^3 \chi_2^3\lambda \in \If_7$. 
\qed \end{coro}

\subsection{The family $\If_{10}$} Recall that
\begin{align*}
\If_{10} &  = \{ \lambda \in \widehat{\Gamma} \mid \lambda_1\notin \G_{12}, \, \lambda_2=1 \}.
\end{align*}

\begin{lemma}\label{lem:modulo irr caso 10}
If $\lambda \in \If_{10}$, then $\dim L(\lambda)=72$. A basis of $L(\lambda)$ is given by
\begin{align*}
\mathrm{B}_{10}&=\{n_{0,b,c,d,e}| c \in \I_{0,1}, b \in\I_{0,3}, d,e \in \I_{0,2} \}.
\end{align*}
\end{lemma}
\pf
Let $w = \widetilde{n}_{1,0,0,0,0}$ and   
$L'(\lambda)=M(\lambda)/\Uc^+ w$. We identify $\mathrm{B}_{10}$ with a basis of $L'(\lambda)$. 
Now $F_1^2 F_{112}^2F_{11212}F_{12}^3  n_{0,3,1,2,2}\in \ku^{\times}v_\lambda$, hence $L'(\lambda)$ is simple.
\epf

Exactly as for Corollary \ref{cor:modulo irr caso 2}, we conclude:
\begin{coro}\label{cor:modulo irr caso 10}
	If $\lambda \in \If_{10}$, then $N(\lambda) \simeq L(\chi_2\lambda)$ and $\chi_2\lambda \in \If_{10}$. 
	\qed
\end{coro}

\section{Simple $\Uc$-modules in class $\cC_2$}\label{sec:irreps2}

We start by the method to compute $L(\lambda)$, $\lambda\in\If_j$, $j\in\I_{11,47}$. We illustrate 
by considering $\If_{11}$, which is small enough to allow complete details; and $\If_{13}$,
with less explicit yet complete enough arguments.
Then we give the main features of the proofs for the other families in $\cC_2$.
Here are the steps of the method:

\begin{enumerate}[leftmargin=*]
\smallbreak\item We identify easily a proper submodule $W=\Uc w_1$ of $M(\lambda)$ as follows:

\begin{itemize} [leftmargin=*]\renewcommand{\labelitemi}{$\diamond$}
\smallbreak\item if $j\in\I_{11,17}$, then $w_1=  \widetilde{m}_{0,0,0,0,1}$, so $W=W_1(\lambda)$, see Lemma \ref{coro:Fi actuando en Wj};

\smallbreak\item if $j\in\I_{18,24}$, then $w_1= \widetilde{m}_{0,0,0,0,2}$, so $W=W_2(\lambda)$, again by Lemma \ref{coro:Fi actuando en Wj};

\smallbreak\item if $j\in\I_{25,35}$, then $w_1$ is as in one of the Remarks \ref{rem:rel E112}, \ref{rem:rel E112 cuadrado},  \ref{rem:rel E12}, \ref{rem:rel E12 cuadrado}, \ref{rem:rel E12 cubo};

\smallbreak\item if $j\in\I_{36,47}$, then $w_1= \widetilde{n}_{1,0,0,0,0}$, so $W=W(\lambda)$ by Lemma \ref{coro:Fi actuando en Wj}.
\end{itemize}

\medbreak
A basis of $M(\lambda)/W$ is obtained by restriction of the height of a specific PBW generator.
Below we denote by $w_2$ an element of $M(\lambda)$ or its class modulo $W$, indistinctly.

\smallbreak\item Next we show that there exists $\beta\neq 0$ and $w_2\in (M(\lambda)/W)_\beta - 0$, such that $F_i w_2=0$, $i\in\I_2$; for this, 
we either apply  one of Remarks \ref{rem:rel E112}, \ref{rem:rel E112 cuadrado}, \ref{rem:l1 a la 3 l2 a la 2}, \ref{rem:rel E12}, \ref{rem:rel E12 cuadrado}, \ref{rem:rel E12 cubo}, or else proceed by direct computation. Hence $\Uc w_2$ is a proper submodule of $M(\lambda)/W$.

\smallbreak\item Let $L'(\lambda)= M(\lambda)/(W + \Uc w_2)$. We consider a suitable set $B_j$ inside the image of the PBW-basis in $L'(\lambda)$ that spans
$L'(\lambda)$. To prove that $\mathrm{B}_j$ is linearly independent, we apply one of the following procedures:

\begin{enumerate}[leftmargin=*, label=\rm{(\alph*)}]
\smallbreak\item\label{it:Bi-li-case1} For $j\in \mathbb{J} =\{ 11,12,18,38\}$, the elements of $\mathrm{B}_j$ are homogeneous of different degrees.

\smallbreak\item\label{it:Bi-li-case2} Assume that $j\notin \mathbb{J}$. Then  $\Uc w_2 \leq M(\lambda)/W$ projects onto the simple module $L(\nu)$,
where $\nu$ is the weight of $w_2$. Also,
let $u \in M(\lambda)/W$ be the element of maximal degree; then $(\Uc u)^{\varphi}$ projects onto a simple  $L(\mu)$.  
Let $\If_k$ and $\If_{\ell}$ be the families containing $\nu$ and $\mu$, respectively. 
At this point we observe that we are proceeding recursively, 
so that we already know the simple modules in $\If_k$ and $\If_{\ell}$.
With this information at hand, we check that $\Uc u = \Uc w_2 \simeq L(\nu)$. 
This isomorphism  provides a basis of $\Uc w_2$; we conclude that there is a  
linear complement of $\Uc w_2$ with a basis $\widetilde{\mathrm{B}}_j$ 
 projecting onto $\mathrm{B}_j$; thus $\mathrm{B}_j$ is a basis of $L'(\lambda)$.
\end{enumerate}

\smallbreak\item Finally we prove that $L'(\lambda)$ is simple. Let $v$ be the element of maximal degree of $L'(\lambda)$. 
A short computation shows that $v$ belongs to every submodule of $L'(\lambda)$.
Applying Lemma \ref{coro:modulos uqsl2} (or by direct computation when we have a table for the action), there exists $F\in\Uc$ such that $Fv=v_\lambda$. 
Hence $L'(\lambda)$ is simple.
\end{enumerate}

As said, we proceed recursively, but with respect to an ad-hoc partial ordering of the families in $\cC_2$. 
In the quiver below, we describe this ordering; $\xymatrix@R=5pt{\If_{11}\ar@{->}[r] & \If_{16}}$
means that knowledge on $\If_{11}$ is used for $\If_{16}$. As we see, there is no vicious circle.

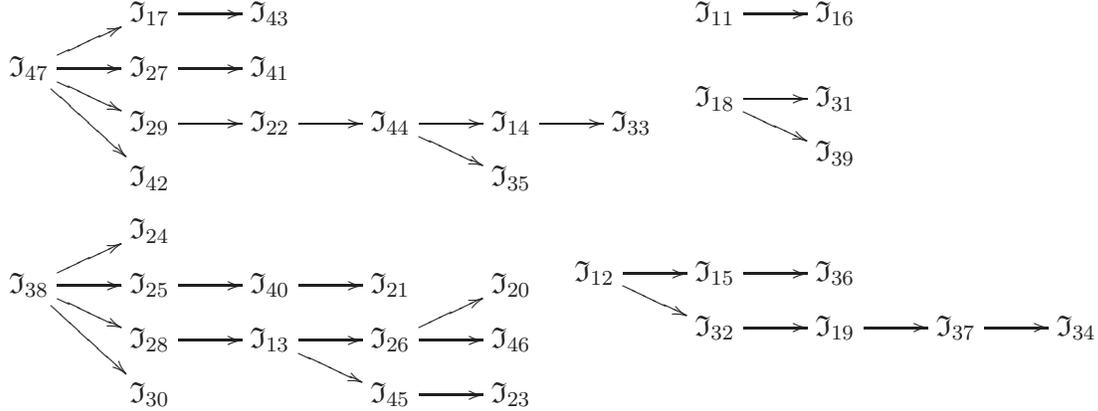
\begin{align*}
&\begin{aligned}
& \xymatrix@R=5pt{
& \If_{17}\ar@{->}[r] & \If_{43} & & & \\
\If_{47} \ar@{->}[r] \ar@{->}[ru] \ar@{->}[rd] \ar@{->}[rdd] & \If_{27}\ar@{->}[r] & \If_{41} & & & \\
& \If_{29}\ar@{->}[r] & \If_{22} \ar@{->}[r] & \If_{44} \ar@{->}[r] \ar@{->}[rd] & \If_{14} \ar@{->}[r] & \If_{33} \\
& \If_{42} & & & \If_{35} &  }
&
& \xymatrix@R=5pt{
\If_{11}\ar@{->}[r] & \If_{16}
\\
\\
\If_{18} \ar@{->}[r] \ar@{->}[rd] & \If_{31} \\
& \If_{39} }
\end{aligned}
\\
& \begin{aligned}
&\xymatrix@R=5pt{
& \If_{24} & & &\\
\If_{38} \ar@{->}[r] \ar@{->}[ru] \ar@{->}[rd] \ar@{->}[rdd] & \If_{25}\ar@{->}[r] & \If_{40} \ar@{->}[r] & \If_{21} & \If_{20}\\
& \If_{28}\ar@{->}[r] & \If_{13} \ar@{->}[r] \ar@{->}[rd] & \If_{26} \ar@{->}[r] \ar@{->}[ru] & \If_{46} \\
& \If_{30} & & \If_{45} \ar@{->}[r]& \If_{23}}
&
& \xymatrix@R=5pt{ \\
\If_{12} \ar@{->}[r] \ar@{->}[rd] & \If_{15}\ar@{->}[r] & \If_{36} & & \\
& \If_{32}\ar@{->}[r] & \If_{19} \ar@{->}[r] & \If_{37} \ar@{->}[r]& \If_{34}}
\end{aligned}
\end{align*}

%In some of the next Lemmas we will describe explicitly the action of $\Uc$ on a fixed basis, whose elements $v_{i,j}$ we define accordingly; these are homogeneous of degree $i\alpha_1+j\alpha_2$.

\subsection{The family $\If_{11}$}
Recall that $\If_{11}   = \{ \lambda \in \widehat{\Gamma} \mid \lambda_1 = 1, \, \lambda_2=\zeta \}$.

\begin{lemma}\label{lem:modulo irr caso 11}
If $\lambda \in \If_{11}$, then $\dim L(\lambda)=11$. A basis of $L(\lambda)$ is given by
\begin{align*}
\mathrm{B}_{11}&=\{m_{a,b,0,d,0}| a \in \I_{0,1}, b \in\I_{0,1}, d \in \I_{0,2} \}-\{m_{1,1,0,0,0}\}.
\end{align*}
The action of $E_i$, $F_i$, $i\in \I_2$ is described in Table \ref{tab:caso 11}.
\end{lemma}

\begin{table}[h]
	\caption{Simple modules for $\lambda \in \If_{11}$}\label{tab:caso 11}
	\begin{center}
		\begin{tabular}{| c | c | c | c | c |} \hline
			$w$ & $E_1\cdot w$ & $E_2\cdot w$ & $\lambda(g_1^{-1})F_1\cdot w$ & $\lambda(g_2^{-1})F_2\cdot w$ \\ \hline
			$v_{0,0}$ & $0$ & $v_{0,1}$ & $0$ & $0$ \\ \hline
			$v_{0,1}$ & $v_{1,1}$ & $0$ & $0$ & $(\zeta^{11}-1)v_{0,0}$ \\ \hline
			$v_{1,1}$ & $v_{2,1}$ & $0$ & $q_{12}(\zeta-1)v_{0,1}$ & $0$\\ \hline
			$v_{2,1}$ & $0$ & $v_{2,2}$ & $q_{12}\zeta^8(1+\zeta^3)v_{1,1}$ & $0$ \\ \hline
			$v_{2,2}$ & $v_{3,2}$ & $0$ & $0$ & $q_{21}^2(1-\zeta)v_{2,1}$\\  \hline
			$v_{3,2}$ & $v_ {4,2}$ & $v_{3,3}$ & $q_{12}^2(\zeta^2-1)v_{2,2}$ & $0$\\ \hline
			$v_{4,2}$ & $0$ & $v_{4,3}$ & $2q_{12}^2(\zeta^2-1)v_{3,2}$ & $0$\\ \hline
			$v_{3,3}$ & $q_{12}\frac{\zeta^8(\zeta^3-1)}{2}v_{4,3}$ & $0$ & $0$ & $q_{21}^3(\zeta^2-1)v_{3,2}$\\  \hline
			$v_{4,3}$ & $v_{5,3}$ & $0$ & $2q_{12}^2(\zeta^2-1)v_{3,3}$ & $q_{21}^4(\zeta^3-1)v_{4,2}$\\ \hline
			$v_{5,3}$ & $0$ & $v_{5,4}$ & $q_{12}^3\zeta^8(1-\zeta^{11})v_{4,3}$ & $0$\\ \hline
			$v_{5,4}$ & $0$ & $0$ & $0$ & $q_{21}^5(\zeta^{11}+1)v_{5,3}$\\ \hline
		\end{tabular}
	\end{center}
\end{table}

\pf Let $w_1=\widetilde{m}_{0,0,0,0,1}$, $w_2=\widetilde{m}_{1,1,0,0,0}$; hence 
$F_i w_1=0$, $i\in\I_2$,
\begin{align*}
F_1  \widetilde{m}_{1,1,0,0,0}& =0, & F_2  \widetilde{m}_{1,1,0,0,0}& = (\zeta^{11}-1)\lambda(g_2) \widetilde{m}_{1,0,0,0,1}\in W_1(\lambda)=\Uc w_1.
\end{align*}
Thus $\Uc w_1+\Uc w_2$ is a proper submodule.
We claim that $ L'(\lambda)=M(\lambda)/\Uc w_1+\Uc w_2$ is simple. 
Consider the following elements of $L'(\lambda)$:
\begin{align*}
v_{0,0}&=\widetilde{m}_{0,0,0,0,0}, & v_{0,1}&=\widetilde{m}_{1,0,0,0,0}, & v_{1,1}&=\widetilde{m}_{0,1,0,0,0}, & v_{2,1}&=\widetilde{m}_{0,0,0,1,0}, \\
v_{2,2}&=\widetilde{m}_{1,0,0,1,0}, & v_{3,2}&=\widetilde{m}_{0,1,0,1,0}, & v_{4,2}&=\widetilde{m}_{0,0,0,2,0}, & v_{3,3}&=\widetilde{m}_{1,1,0,1,0}, \\
v_{4,3}&=\widetilde{m}_{1,0,0,2,0}, & v_{5,3}&=\widetilde{m}_{0,1,0,2,0}, & v_{5,4}&=\widetilde{m}_{1,1,0,2,0}.
\end{align*}
Notice that $v_{i,j} \in L'(\lambda)_{i\alpha_1+j\alpha_2}$. The action of $E_i$, $F_i$ on these vectors is given in Table \ref{tab:caso 11}, and we check that $L'(\lambda)$ is spanned by the $v_{i,j}$'s by direct computation. 

%Note that $E_1  v_{5,4}=E_2  v_{5,4}=0$. 
For each $v_{i,j}$ there exist $E_{i,j} \in\Uc^+_{(5-i)\alpha_1+(4-j)\alpha_2}$ such that $E_{i,j}  v_{i,j}=v_{5,4}$; also, there exists $F_{5,4}\in\Uc^-_{-5\alpha_1-4\alpha_2}$ such that $F_{5,4}  v_{5,4}=v_\lambda$. This implies that the $v_{i,j}$'s are $\neq 0$; hence they are linearly independent, since they have different degrees, and $\mathrm{B}_{11}$ is identified with a basis of $L'(\lambda)$.

Let now $0\neq U\leq L'(\lambda)$ and pick $v\in U - 0$. Expressing $v$ in the basis $\mathrm{B}_{11}$, we see that 
there exists $E\in \Uc^+$ such that $E  v=v_{5,4}$. But $\Uc  v_{5,4}=L'(\lambda)$. Hence $L'(\lambda)$ is simple.
\epf

\begin{remark}\label{cor:modulo irr caso 11}
If $\lambda \in \If_{11}$, then $N(\lambda)/W_1(\lambda) \simeq L(\chi_1 \chi_2^2 \lambda)$, with $\chi_1 \chi_2^2 \lambda \in \If_{41}$ has dimension 37. 
Now $W_1(\lambda)$ is a lowest weight module of lowest weight $ \chi_1 \lambda \in \If_{43}$; since $\dim L( \chi_1 \lambda) = 25$ by Lemma \ref{lem:modulo irr caso 43},
the kernel of $W_1(\lambda) \twoheadrightarrow L( \chi_1 \lambda)$ is a submodule of dimension 71.
\end{remark}

\subsection{The family $\If_{12}$} Recall that $\If_{12}   = \{ \lambda \in \widehat{\Gamma} \mid \lambda_1 = 1, \, \lambda_2=\zeta^4 \}$.

\begin{lemma}\label{lem:modulo irr caso 12}
If $\lambda \in \If_{12}$, then $\dim L(\lambda)=11$. A basis of $L(\lambda)$ is given by
\begin{align*}
\mathrm{B}_{12}&=\{m_{a,b,0,d,0}: a,b,d \in \I_{0,1}\}\cup \{m_{0,1,1,0,0},m_{1,0,1,1,0},m_{0,0,1,1,0}\}. 
\end{align*}
The action of $E_i$, $F_i$, $i\in \I_2$ is described in Table \ref{tab:caso 12}.
\end{lemma}

\begin{table}[h]
	\caption{Simple modules for $\lambda \in \If_{12}$}\label{tab:caso 12}
	\begin{center}
		\begin{tabular}{| c | c | c | c | c |} \hline
			$w$ & $E_1\cdot w$ & $E_2\cdot w$ & $\lambda(g_1^{-1})F_1\cdot w$ & $\lambda(g_2^{-1})F_2\cdot w$ \\ \hline
			$v_{0,0}$ & $0$ & $v_{0,1}$ & $0$ & $0$ \\ \hline
			$v_{0,1}$ & $v_{1,1}$ & $0$ & $0$ & $(\zeta^{10}+1)v_{0,0}$ \\ \hline
			$v_{1,1}$ & $v_{2,1}$ & $v_{1,2}$ & $q_{12}(\zeta-1)v_{0,1}$ & $0$\\ \hline
			$v_{2,1}$ & $0$ & $v_{2,2}$ & $q_{12}\zeta^8(1+\zeta^3)v_{1,1}$ & $0$ \\ \hline
			$v_{1,2}$ & $\zeta^{11}(1+\zeta^3)q_{12}v_{2,2}$ & $0$ & $0$ & $q_{21}(1+\zeta^3)\zeta^4v_{1,1}$\\  \hline
			$v_{2,2}$ & $v_ {3,2}$ & $0$ & $q_{12}(\zeta^3+1)\zeta^8v_{1,2}$ & $-q_{2,1}^2v_{2,1}$\\ \hline
			$v_{3,2}$ & $0$ & $v_{3,3}$ & $q_{12}^2\zeta^{10}v_{2,2}$ & $0$\\ \hline
			$v_{3,3}$ & $0$ & $0$ & $0$ &  $q_{21}^3\zeta^3(1-\zeta)v_{3,2}$ \\  \hline
			$v_{4,3}$ & $\zeta^9q_{12}v_{5,3}$ & $0$ & $q_{12}^4\zeta(3)_{\zeta^{11}} v_{3,3}$ & $0$\\ \hline
			$v_{5,3}$ & $0$ & $v_{5,4}$ & $-q_{12}^2(1+\zeta^{3})v_{4,3}$ & $0$\\ \hline
			$v_{5,4}$ & $0$ & $0$ & $0$ & $q_{21}^5(1-\zeta)\zeta^4v_{5,3}$\\ \hline
		\end{tabular}
	\end{center}
\end{table}

\pf 
Let $w_1=\widetilde{m}_{0,0,0,0,1}$, $w_2=F_2E_2E_{12}^2v_\lambda$;
then $F_iw_j=0$ for $i,j\in\I_2$, so $\Uc w+ W_1(\lambda)$ is a proper submodule of $M(\lambda)$.
Let $L'(\lambda) := M(\lambda)/\Uc w+ W_1(\lambda)$. We label the elements of $\mathrm{B}_{12}$ as follows:
\begin{align*}
v_{0,0}&=m_{0,0,0,0,0}, & v_{0,1}&=m_{1,0,0,0,0}, & v_{1,1}&=m_{0,1,0,0,0}, & v_{2,1}&=m_{0,0,0,1,0}, \\
v_{2,2}&=m_{1,0,0,1,0}, & v_{1,2}&=m_{1,1,0,0,0},& v_{3,2}&=m_{0,1,0,1,0},  & v_{3,3}&=m_{1,1,0,1,0}, \\
v_{4,3}&=m_{0,1,1,0,0}, & v_{5,3}&=m_{0,0,1,1,0}, & v_{5,4}&=m_{1,0,1,1,0}.
\end{align*}
The action of $E_i$, $F_i$ on these vectors is given in Table and $\mathrm{B}_{12}$ is a basis of $L'(\lambda)$. 
Looking at the table, there exists $F\in\Uc^-$ such that $F m_{1,0,1,1,0}=v_\lambda$.
Then $L'(\lambda)$ is simple.
\epf

\subsection{The family $\If_{13}$} Recall that $\If_{13}   = \{ \lambda \in \widehat{\Gamma} \mid \lambda_1 = 1, \, \lambda_2=\zeta^7 \}$.

\begin{lemma}\label{lem:modulo irr caso 13}
If $\lambda \in \If_{13}$, then $\dim L(\lambda)=23$. A basis of $L(\lambda)$ is given by
\begin{align*}
\mathrm{B}_{13}&=\{m_{a,b,0,d,0}|b\in \I_{0,2}\}\cup\{m_{a,0,1,0,0},m_{0,3,0,d,0},m_{1,3,0,1,0}|a\in \I_{0,1},   d \in \I_{1,2}\}. 
\end{align*}
\end{lemma}
\pf
Let $w_1=\widetilde{m}_{0,0,0,0,1}$, $w_2 =F_2 E_2 E_{12}^3v_\lambda$. Then 
$W_1(\lambda)=\Uc w_1$ by Lemma \ref{coro:Fi actuando en Wj}, and $F_1w_2=F_2w_2=0$ by Remark 4.22, so $\Uc w_1+\Uc w_2 \lneq M(\lambda)$. We claim that 
$L'(\lambda):= M(\lambda)/(\Uc w_1+ \Uc w_2)$ is simple and $\mathrm{B}_{13}$ is a basis of $L'(\lambda)$.

Let $M=M(\lambda)/W_1(\lambda)$ and $u=m_{1,3,1,2,0} \in M$. Notice that
$E_{112}^2E_{11212}E_2 w_2 = -q_{12}^{18} \ u$, so $u\in\Uc w_2$.
On the other hand, $E_iu=0$, $i\in \I_2$, $g_1\gpl_1u=u$ and $g_2\gpl_2u=\zeta^9u$,
so $(\Uc u)^\varphi$ projects over a simple module $L(\mu)$ with $\mu\in \If_{14}$, see Lemma \ref{lem: irreduble via phi}; in particular there exists
$F'\in\Uc_{-7\alpha_1-5\alpha_2}$ such that $F'u\neq 0$. As $\Uc u\subseteq \Uc w_2$ and $\Uc w_2$ is a lowest weight module,
\begin{align*}
F'u &\in (\Uc u)_{3\alpha_1+3\alpha_2}\subseteq (\Uc w_2)_{3\alpha_1+3\alpha_2}=\ku w.
\end{align*}
Hence we may assume that $F'u=w_2$, and $\Uc u=\Uc w_2$.

Also $g_1\gpl_1w_2=\zeta^9w_2$, $g_2\gpl_2w_2 =\zeta^4w_2$, so $\Uc w_2$ projects over a simple module $L(\nu)$ with $\nu\in \If_{28}$. 
For any $v\in M$, $v\neq0$, there exists $E\in\Uc$ such that $Ev=u$. Thus we conclude that $\Uc w_2 \simeq L(\nu)$, and then $\dim L'(\lambda)=48-25=23$
by Lemma \ref{lem:modulo irr caso 28}.

%From $w$ and applying repeatedly $E_1,E_2$ over $w$ we obtain
%\begin{align*}
%m_{0,3,0,0,0}& =\frac{q_{12}^2\zeta^5(4)_{\zeta^7}}{3} m_{1,1,0,1,0}+ \frac{q_{12}\zeta(1+\zeta^2)(3)_{\zeta^5}}{3} m_{1,0,1,0,0},  \\
%m_{0,1,1,0,0}&=q_{12}(1+\zeta^3)m_{0,2,0,1,0}+q_{12}^2m_{1,0,0,2,0},\\ 
%m_{0,0,1,1,0}&=q_{12}\zeta^7(1+\zeta^2)m_{0,1,0,2,0},\\
%m_{0,2,1,0,0}&=q_{12}^3\zeta^{10}(1-\zeta) m_{1,1,0,2,0},\\ 
%m_{0,3,1,0,0}&=\zeta^4q_{12}^4(1-\zeta)m_{1,2,0,2,0},\\
%m_{0,2,1,1,0}&=q_{12}\zeta^7(1+\zeta^2)m_{0,3,0,2,0}. 
%\end{align*}

Applying Lemma \ref{coro:modulos uqsl2}, there exists $F\in\Uc^-$ such that $F  m_{0,3,0,2,0}=v_\lambda$. Note that $$E_2 m_{0,3,0,2,0}= m_{1,3,0,2,0}=0$$ since $0=E_{12} m_{0,3,1,0,0}$ and $\ku m_{1,2,1,1,0}=\ku m_{1,3,0,2,0}$. Also $E_1 m_{0,3,0,2,0} =0 $ because it is a scalar multiple of $m_{0,1,1,2,0}$ , which is $0$. Using this fact and previous relations, we are able to prove that $\mathrm{B}_{13}$ spans $L'(\lambda)$, but as $\mathrm{B}_{13}$ has 23 elements, it is a basis.

Let $0\neq W \leq L'(\lambda)$, $w\in W-0$. Arguing as before, there exists $E\in\Uc^+$ such that $E  w=m_{0,3,0,2,0}$, so $m_{0,3,0,2,0}\in W$, but then $v_\lambda\in W$, so  $L'(\lambda)$ is simple.
\epf

\subsection{The family $\If_{14}$} Recall that $\If_{14}   = \{ \lambda \in \widehat{\Gamma} \mid \lambda_1 = 1, \, \lambda_2=\zeta^3 \}$.

\begin{lemma}\label{lem:modulo irr caso 14}
If $\lambda \in \If_{14}$, then $\dim L(\lambda)=25$. A basis of $L(\lambda)$ is given by
\begin{align*}
\mathrm{B}_{14}&=\{m_{a,b,0,d,0}\mid a \in \I_{0,1}, b \in\I_{0,3}, d \in \I_{0,2} \} \cup \{m_{0,0,1,0,0},m_{0,0,1,2,0}\}-\{m_{1,3,0,2,0}\}. 
\end{align*}
\end{lemma}

\pf
Let $w_1=\widetilde{m}_{0,0,0,0,1}$, $w_2 =(1+\zeta^3) \widetilde{m}_{1,0,1,0,0}+ q_{12}\zeta^3(1+\zeta)\widetilde{m}_{1,1,0,1,0}$. Then $W_1(\lambda)=\Uc w_1$ and $F_1w_2=F_2w_2=0$ by direct computation. 

Let $M=M(\lambda)/W_1(\lambda)$, $L'(\lambda)=M(\lambda)/\Uc w_2+ W_1(\lambda)$ and $u=m_{1,3,1,2,0} \in M$. 
Then $(\Uc u)^\varphi$ projects over $L(\mu)$ for some $\mu\in \If_{13}$. 
Also, $\Uc w_2$ projects over $L(\nu)$ for some $\nu\in \If_{44}$. Hence $\Uc u=\Uc w_2$, and moreover $\Uc w_2$ is simple, so $\dim L'(\lambda)=48-25=23$
by Lemma \ref{lem:modulo irr caso 44}.
By direct computation $L'(\lambda)$ is spanned by $\mathrm{B}_{14}$, so $\mathrm{B}_{14}$ is a basis of $L'(\lambda)$. 

Moreover there exists $F\in\Uc^-$ such that $F m_{1,0,1,2,0}=v_\lambda$, so $L'(\lambda)$ is simple.
\epf

\subsection{The family $\If_{15}$} Recall that $\If_{15}   = \{ \lambda \in \widehat{\Gamma} \mid \lambda_1 = 1, \, \lambda_2=\zeta^9 \}$.

\begin{lemma}\label{lem:modulo irr caso 15}
If $\lambda \in \If_{15}$, then $\dim L(\lambda)=37$. A basis of $L(\lambda)$ is given by
\begin{align*}
\mathrm{B}_{15}&=\{m_{a,b,c,d,0}| a,c \in \I_{0,1}, b \in\I_{0,3}, d \in \I_{0,2}  \}
\\ &-\{m_{a,b,1,d,0}| a \in \I_{0,1}, b  \in \I_{2,3}, d \in \I_{0,2}, (a,b,d)\neq(0,2,2)\}. 
\end{align*}
\end{lemma}

\pf
Let $w_1=\widetilde{m}_{0,0,0,0,1}$, $u=\widetilde{m}_{1,3,1,2,0}$, $w_2 =F_2F_{12}F_{112}^2u$. Then 
$W_1(\lambda)=\Uc w_1$.

Let $M=M(\lambda)/W_1(\lambda)$, so $E_1u=E_2u=0$ in $M$, and $(\mathcal U u)^{\varphi}\twoheadrightarrow L(\nu)$  
for some  $\nu\in \If_{11}$; thus $w_2\neq 0$. By direct computation, $F_iw_2=0$, $i\in\I_2$, so $\mathcal U w_2$ projects over a simple module $L(\mu)$, for $\mu \in \If_{12}$. From here, $\mathcal U w_2 \simeq L(\mu)$.

Let $L'(\lambda)=M(\lambda)/W_1(\lambda)+\Uc w_2$. Then $\dim L'(\lambda)=37$ by Lemma \ref{lem:modulo irr caso 12}, 
and $\mathrm{B}_{15}$ is a basis of $L'(\lambda)$. There exists $F$ such that $Fm_{0,2,1,2,0}=v_\lambda$, and $L'(\lambda)$ is simple.
\epf

\subsection{The family $\If_{16}$} Recall that $\If_{16}   
= \{ \lambda \in \widehat{\Gamma} \mid \lambda_1 = 1, \, \lambda_2= -1 \}$.

\begin{lemma}\label{lem:modulo irr caso 16}
If $\lambda \in \If_{16}$, then $\dim L(\lambda)=37$. A basis of $L(\lambda)$ is given by
\begin{align*}
\mathrm{B}_{16}&=\{m_{a,b,c,d,0} | a,c \in \I_{0,1}, b \in\I_{0,3}, d \in \I_{0,2} \}
\\ &- 
\big( \{ m_{a,3,c,d,0}|  a,c \in \I_{0,1}, d \in \I_{1,2}  \}   \cup  \{m_{1,2,1,2,0}, m_{0,2,1,2,0}, m_{1,2,0,2,0} \} \big).
\end{align*}
\end{lemma}

\pf
Let $w_1=\widetilde{m}_{0,0,0,0,1}$, $u=\widetilde{m}_{1,3,1,2,0}$, $w_2 =F_2F_{11212}F_{112}u$. Then 
$W_1(\lambda)=\Uc w_1$.

Let $M=M(\lambda)/W_1(\lambda)$, so $E_1u=E_2u=0$ in $M'$, and $(\mathcal U u)^{\varphi}\twoheadrightarrow L(\nu)$ for some  $\nu\in \If_{12}$; thus $w_2\neq 0$. By direct computation, $F_iw_2=0$, $i\in\I_2$, so $\mathcal U w_2$ projects over a simple module $L(\mu)$, for $\mu \in \If_{11}$. From here, $\mathcal U w_2\simeq L(\mu)$.

Let $L'(\lambda)=M(\lambda)/ W_1(\lambda) + \Uc w_2$.
Then $\dim L'(\lambda)=37$ by Lemma \ref{lem:modulo irr caso 11}, and $\mathrm{B}_{16}$ is a basis of $L'(\lambda)$. There exists $F$ such that $Fm_{1,1,1,2,0}=v_\lambda$, and $L'(\lambda)$ is simple.
\epf

\subsection{The family $\If_{17}$} Recall that $\If_{17}   = \{ \lambda \in \widehat{\Gamma} \mid \lambda_1 = 1, \, \lambda_2=\zeta^{10} \}$.

\begin{lemma}\label{lem:modulo irr caso 17}
If $\lambda \in \If_{17}$, then $\dim L(\lambda)=47$. A basis of $L(\lambda)$ is given by
\begin{align*}
\mathrm{B}_{17}&=\{m_{a,b,c,d,0}|  a,c \in \I_{0,1}, b \in\I_{0,3}, d,e \in \I_{0,2}, 	(a,b,c,d) \neq (1,3,1,2) \}. 
\end{align*}
\end{lemma}

\pf
Let $w_1=\widetilde{m}_{0,0,0,0,1}$, $w_2 =\widetilde{m}_{1,3,1,2,0}$. Then 
$W_1(\lambda)=\Uc w_1$, and $F_iw=0$, $i\in\I_2$, so $\mathcal U w$ projects over a simple module $L(\mu)$, for $\mu \in \If_{47}$.
Let $M=M(\lambda)/W_1(\lambda)$, hence $\mathcal U w_2\simeq L(\mu)$.
Let $L'(\lambda)=M(\lambda)/ W_1(\lambda) + \Uc w_2$.
Then $\dim L'(\lambda)=47$ by Lemma \ref{lem:modulo irr caso 47}, and $\mathrm{B}_{17}$ is a basis of $L'(\lambda)$. There exists $F$ such that $Fm_{0,3,1,2,0}=v_\lambda$, and $L'(\lambda)$ is simple.
\epf

\subsection{The family $\If_{18}$} Recall that $\If_{18}   = \{ \lambda \in \widehat{\Gamma} \mid \lambda_1 = \zeta^8, \, \lambda_2=\zeta^5 \}$.

\begin{lemma}\label{lem:modulo irr caso 18}
If $\lambda \in \If_{18}$, then $\dim L(\lambda)=11$. A basis of $L(\lambda)$ is given by  \begin{align*}
\mathrm{B}_{18}&=\{m_{a,b,1,0,1}|a,b \in \I_{0,1}\} \cup \{m_{0,b,0,0,e}|e \in \I_{0,1}, b \in\I_{0,3}\}
\cup \{m_{1,0,0,0,0}\}
\\ &-\{m_{1,1,1,0,1}, m_{3,0,0,0,1}\}. 
\end{align*}
The action of $E_i$, $F_i$, $i\in \I_2$ is described in Table \ref{tab:caso 18}.
\end{lemma}

\begin{table}[h]
	\caption{Simple modules for $\lambda \in \If_{18}$}\label{tab:caso 18}
	\begin{center}
		\begin{tabular}{| c | c | c | c | c |} \hline
			$w$ & $E_1\cdot w$ & $E_2\cdot w$ & $\lambda(\gpl_1)F_1\cdot w$ & $\lambda(g_2)^{-1}F_2\cdot w$ \\ \hline
			$v_{0,0}$ & $v_{1,0}$ & $v_{0,1}$ & $0$ & $0$ \\ \hline
			$v_{1,0}$ & $0$ & $q_{21}\zeta^9(4)_{\zeta} v_{1,1}$ & $(1+\zeta^2)v_{0,0}$ & $0$ \\ \hline
			$v_{0,1}$ & $\zeta^8(4)_{\zeta}v_{1,1}$ & 0 & $0$ & $(\zeta^7-1)v_{0,0}$ \\ \hline
			$v_{1,1}$ & $\frac{q_{12}\zeta^4(4)_{\zeta^7}}{3}v_{2,1}$ & $0$ & $q_{12}(\zeta-1)v_{0,1}$ & $(\zeta^{11}-1)v_{1,0}$\\ \hline
			$v_{2,1}$ & $0$ & $q_{21}^2\zeta^{10}(4)_{\zeta}v_{2,2}$ & $(1-\zeta^4) v_{1,1}$ & $0$ \\ \hline
			$v_{2,2}$ & $(1-\zeta^4)v_{3,2}$ & $0$ & $0$ & $\frac{-(1+\zeta^2)(3)_{\zeta^7}}{3}v_{2,1}$\\ \hline
			$v_{3,2}$ & $v_{4,2}$ & $q_{12}\zeta^{10}(4)_{\zeta}v_{3,3}$ & $\zeta^{10}(4)_{\zeta}v_{2,2}$ & $0$\\\hline
			$v_{4,2}$ & $0$ & $v_{4,3}$ & $q_{12}^2\zeta(\zeta+1)v_{3,2}$ & $0$\\\hline
			$v_{3,3}$ & $\frac{q_{12}^4\zeta^7(4)_{\zeta}}{3}v_{4,3}$ & $0$ & $0$ & $\frac{\zeta^8-1}{3}v_{3,2}$\\\hline
			$v_{4,3}$ & $v_{5,3}$ & $0$ & $q_{12}^3(\zeta^{11}+1)(4)_{\zeta}^2 v_{3,3}$ & $q_{21}^4(\zeta^{11}-1)v_{4,2}$\\\hline
			$v_{5,3}$ & $0$ & $0$ & $q_{12}^3\zeta^4v_{4,3}$ & $0$\\\hline
		\end{tabular}
	\end{center}
\end{table}

\pf $W_2(\lambda)\leq M(\lambda)$ by Lemma \ref{coro:Fi actuando en Wj} and $w:=F_2E_2E_{12} $ satisfies $F_1w=F_2w=0$ by Remark \ref{rem:rel E12}. Let $L'(\lambda)=M(\lambda)/\Uc w_2 +W_2(\lambda)$. We fix the following notation for  $\mathrm{B}_{18}$:
\begin{align*}
v_{0,0}&=m_{0,0,0,0,0}, & v_{1,0}&=m_{0,0,0,0,1},& v_{0,1}&=m_{1,0,0,0,0}, & v_{1,1}&=m_{0,1,0,0,0}, \\
v_{2,1}&=m_{0,1,0,0,1}, & v_{2,2}&=m_{0,2,0,0,0}, & v_{3,2}&=m_{0,2,0,0,1}, & v_{4,2}&=m_{0,0,1,0,1}, \\
v_{3,3}&=m_{0,3,0,0,0}, & v_{4,3}&=m_{1,0,1,0,1}, & v_{5,3}&=m_{0,1,1,0,1}.
\end{align*}

We check that $L'(\lambda)$ is spanned by
$\mathrm{B}_{18}$. From Table \ref{tab:caso 18} there exist $E_{i,j} \in\Uc^+_{(5-i)\alpha_1+(3-j)\alpha_2}$,
$F_{5,3}\in\Uc^-_{-5\alpha_1-3\alpha_2}$ such that $E_{i,j}  v_{i,j}=v_{5,3}$, $F_{5,3} v_{5,3}=v_\lambda$.
Thus $L'(\lambda)$ is simple.
\epf

\subsection{The family $\If_{19}$}
Recall that $\If_{19}   = \{ \lambda \in \widehat{\Gamma} \mid \lambda_1 = \zeta^8, \, \lambda_2=\zeta^8 \}$.

\begin{lemma}\label{lem:modulo irr caso 19}
If $\lambda \in \If_{19}$, then $\dim L(\lambda)=35$. A basis of $L(\lambda)$ is given by
\begin{align*}
\mathrm{B}_{19}= &\{m_{0,b,0,d,e}| b \in\I_{0,3}, d\in \I_{0,2}, e \in \I_{0,1} \}  
\cup  \{m_{1,b,0,0,e}| \ b,e  \in \I_{0,1}\} \cup \{m_{0,b,1,0,0}| \ b \in \I_{1,3}\} \\
 & \cup  \{m_{1,b,0,0,1}\mid b \in \I_{2,3}\} \cup \{m_{1,0,0,1,1}, m_{0,0,1,1,0}\}. 
\end{align*}
\end{lemma}

\pf 
Let $w_1=\widetilde{m}_{0,0,0,0,2}$, $w_2 =F_2E_2E_{12}^2v_\lambda$. Then $W_2(\lambda)=\Uc w_1$ and $F_1w_2=F_2w_2=0$.
Set $M'=M(\lambda)/W_2(\lambda)$, $u=\widetilde{m}_{1,3,1,2,1}$. 
Hence $\Uc w_2\twoheadrightarrow L(\mu)$ for $\mu \in \If_{32}$, and there exists $E\in\Uc$ such that
$Ew_2=u$.  Moreover, there exists $F \in \Uc$
such that $Fu=w_2$, so $\Uc w_2 = \Uc u \simeq L(\mu)$.
Let $L'(\lambda)=M(\lambda)/\Uc w_2 + W_2(\lambda)$, so $\dim L'(\lambda)=96-61=35$ by Lemma \ref{lem:modulo irr caso 32}, and $\mathrm{B}_{19}$ is a basis of $L'(\lambda)$. As in previous cases, $L'(\lambda)$ is simple.
\epf

\subsection{The family $\If_{20}$}
Recall that $\If_{20}   = \{ \lambda \in \widehat{\Gamma} \mid \lambda_1 = \zeta^8, \, \lambda_2=\zeta^{11} \}$.

\begin{lemma}\label{lem:modulo irr caso 20}
If $\lambda \in \If_{20}$, then $\dim L(\lambda)=71$. A basis of $L(\lambda)$ is given by
\begin{align*}
\mathrm{B}_{20}= &\{m_{a,b,c,d,e}| a,c, e \in \I_{0,1}, b \in\I_{0,3}, d \in \I_{0,2} \}  
\\&- \Big( \{ m_{1,b,1,d,e} |b \in\I_{0,3}, d \in \I_{0,2}, e \in \I_{0,1}, (b,d,e)\neq(2,2,1) \} 
\cup \{ m_{1,0,0,2,1}, m_{1,3,0,0,0} \} \Big).
\end{align*}
\end{lemma}
\pf
Let $w_1=\widetilde{m}_{0,0,0,0,2}$, $w_2 =F_2E_2E_{12}^3v_\lambda$. Then $W_2(\lambda)=\Uc w_1$ and $F_1w_2=F_2w_2=0$.
Set $M'=M(\lambda)/W_2(\lambda)$, $u=\widetilde{m}_{1,3,1,2,1}$.
Hence $\Uc w_2\twoheadrightarrow L(\mu)$ for $\mu \in \If_{26}$, and there exists $E\in\Uc$ such that
$Ew_2=u$.  Moreover, there exists $F \in \Uc$
such that $Fu=w_2$, so $\Uc w_2=\Uc u \simeq L(\mu)$.
Let $L'(\lambda)=M(\lambda)/\Uc w_2 + W_2(\lambda)$, so $\dim L'(\lambda)=96-25=71$ by Lemma \ref{lem:modulo irr caso 26} and $\mathrm{B}_{20}$ is a basis of $L'(\lambda)$. As in previous cases, $L'(\lambda)$ is simple.
\epf

\subsection{The family $\If_{21}$}
Recall that $\If_{21}   = \{ \lambda \in \widehat{\Gamma} \mid \lambda_1 = \zeta^8, \, \lambda_2=\zeta^3 \}$.

\begin{lemma}\label{lem:modulo irr caso 21}
If $\lambda \in \If_{21}$, then $\dim L(\lambda)=61$. A basis of $L(\lambda)$ is given by
\begin{align*}
\mathrm{B}_{21}= &\{ m_{a,b,c,d,e}\mid  a,b,c, e \in \I_{0,1},  d \in \I_{0,2} \}
 \cup \{m_{a,2,c,0,e}\mid a,c,e \in \I_{0,1}\} \\
 & \cup \{m_{1,3,0,0,e}\mid e \in \I_{0,1} \} \cup \{m_{0,3,1,0,1}, m_{1,3,1,0,1}, m_{0,2,0,1,0}\}.
\end{align*}
\end{lemma}
\pf
Let $w_1=\widetilde{m}_{0,0,0,0,2}$, $u=\widetilde{m}_{1,3,1,2,1}$, $w_2 =F_1F_{11212}F_{12}u$. Then 
$W_2(\lambda)=\Uc w_1$.

Let $M'=M(\lambda)/W_2(\lambda)$, so $E_1u=E_2u=0$ in $M'$, and $(\mathcal U u)^{\varphi}\twoheadrightarrow L(\nu)$ for some  $\nu\in \If_{19}$; thus $w_2\neq 0$.
By direct computation, $F_iw_2=0$, $i\in\I_2$, so $\mathcal{U} w_2$ projects over a simple module $L(\mu)$, for $\mu \in \If_{40}$. From here, $\mathcal{U} w_2\simeq L(\mu)$.

Let $L'(\lambda)=M(\lambda)/ W_2(\lambda) + \Uc w_2$.
Then $\dim L'(\lambda)=61$ by Lemma \ref{lem:modulo irr caso 40}, and $\mathrm{B}_{21}$ is a basis of $L'(\lambda)$. There exists $F$ such that $Fm_{1,1,1,2,1}=v_\lambda$, and $L'(\lambda)$ is simple.
\epf

\subsection{The family $\If_{22}$}
Recall that $\If_{22}   = \{ \lambda \in \widehat{\Gamma} \mid \lambda_1 = \zeta^8, \, \lambda_2=\zeta^9 \}$.

\begin{lemma}\label{lem:modulo irr caso 22}
If $\lambda \in \If_{22}$, then $\dim L(\lambda)=49$. A basis of $L(\lambda)$ is given by
\begin{align*}
\mathrm{B}_{22}= &\{m_{a,b,c,d,e}| \ a,c \in \I_{0,1}, b \in\I_{0,3}, d,e \in \I_{0,1}\} 
\\& - \{ m_{a,b',1,0,0},m_{1,3,1,1,1},m_{a,b,1,1,0} \mid  a \in \I_{0,1}, b' \in\I_{0,3}, b \in\I_{1,3} \}.
\end{align*}
\end{lemma}

\pf
Let $w_1=\widetilde{m}_{0,0,0,0,2}$, $w_2 =F_1^2E_{112}^2E_1 v_\lambda$. Then $W_2(\lambda)=\Uc w_1$ and $F_1w_2=F_2w_2=0$.
Set $M'=M(\lambda)/W_2(\lambda)$, $u=\widetilde{m}_{1,3,1,2,1}$.
Hence $\Uc w_2\twoheadrightarrow L(\mu)$ for $\mu \in \If_{29}$, and there exists $E\in\Uc$ such that
$Ew_2=u$.  Moreover, there exists $F \in \Uc$
such that $Fu=w_2$, so $\Uc w_2=\Uc u \simeq L(\mu)$.
Let $L'(\lambda)=M(\lambda)/\Uc w_2 + W_2(\lambda)$, so $\dim L'(\lambda)=96-47=49$ by Lemma \ref{lem:modulo irr caso 29}, and $\mathrm{B}_{22}$ is a basis of $L'(\lambda)$. As in previous cases, $L'(\lambda)$ is simple.
\epf

\subsection{The family $\If_{23}$}
Recall that $\If_{23}   = \{ \lambda \in \widehat{\Gamma} \mid \lambda_1 = \zeta^8, \, \lambda_2=\zeta^2 \}$.

\begin{lemma}\label{lem:modulo irr caso 23}
If $\lambda \in \If_{23}$, then $\dim L(\lambda)=47$. A basis of $L(\lambda)$ is given by
\begin{align*}
\mathrm{B}_{23}= &\Big(\{m_{a,b,0,d,e}| a, e \in \I_{0,1}, b \in\I_{0,3}, d \in \I_{0,2} \} 
\cup \{m_{a,b,1,0,0} \mid a,b \in \I_{0,1} \} \\
&\cup \{m_{0,2,1,0,0},  m_{1,3,1,0,0}\} \Big)  - \Big(\{ m_{1,b,0,1,e} |b \in\I_{0,2}, e \in \I_{0,1}\} \cup  \{m_{0,2,0,2,0}\} \Big).
\end{align*}
\end{lemma}
\pf
Let $w_1=\widetilde{m}_{0,0,0,0,2}$, $u=\widetilde{m}_{1,3,1,2,1}$, $w_2 =F_{12}^3F_{11212}F_{112}F_1u$. Then 
$W_2(\lambda)=\Uc w_1$.

Let $M'=M(\lambda)/W_2(\lambda)$, so $E_1u=E_2u=0$ in $M'$, and $(\mathcal U u)^{\varphi}\twoheadrightarrow L(\nu)$ for some  $\nu\in \If_{22}$; thus $w_2\neq 0$.
By direct computation, $F_iw_2=0$, $i\in\I_2$, so $\mathcal{U} w_2$ projects over a simple module $L(\mu)$, for $\mu \in \If_{45}$. From here, $\mathcal{U} w_2\simeq L(\mu)$.

Let $L'(\lambda)=M(\lambda)/ W_1(\lambda) + \Uc w_2$.
Then $\dim L'(\lambda)=47$ by Lemma \ref{lem:modulo irr caso 45}, and $\mathrm{B}_{23}$ is a basis of $L'(\lambda)$. There exists $F$ such that $F m_{1,3,0,2,1}=v_\lambda$, and $L'(\lambda)$ is simple.
\epf

\subsection{The family $\If_{24}$}Recall that $\If_{24}   = \{ \lambda \in \widehat{\Gamma} \mid \lambda_1 = \zeta^8, \, \lambda_2= -1 \}$.

\begin{lemma}\label{lem:modulo irr caso 24}
If $\lambda \in \If_{24}$, then $\dim L(\lambda)=85$. A basis of $L(\lambda)$ is given by
\begin{align*}
\mathrm{B}_{24}= &\{m_{a,b,c,d,e}| a,c,e \in \I_{0,1}, b \in\I_{0,3}, d \in \I_{0,2} \} \\ &- \big(\{ m_{a,3,c,2,e},m_{1,3,c,1,1}| a,c,e \in \I_{0,1} \} \cup \{m_{0,3,1,1,1}\}\big).
\end{align*}
\end{lemma}

\pf
Let $w_1=\widetilde{m}_{0,0,0,0,2}$, $u=\widetilde{m}_{1,3,1,2,1}$, $w_2 =F_{12}F_{11212}F_1 u$. Then 
$W_2(\lambda)=\Uc w_1$.

Let $M'=M(\lambda)/W_2(\lambda)$, so $E_1u=E_2u=0$ in $M'$, and $(\mathcal U u)^{\varphi}\twoheadrightarrow L(\nu)$ for some  $\nu\in \If_{18}$; thus $w_2\neq 0$.
By direct computation, $F_iw_2=0$, $i\in\I_2$, so $\mathcal{U} w_2$ projects over a simple module $L(\mu)$, for $\mu \in \If_{38}$. From here, $\mathcal{U} w_2\simeq L(\mu)$.

Let $L'(\lambda)=M(\lambda)/ W_1(\lambda) + \Uc w_2$.
Then $\dim L'(\lambda)=85$ by Lemma \ref{lem:modulo irr caso 38}, and $\mathrm{B}_{24}$ is a basis of $L'(\lambda)$. There exists $F$ such that $F m_{1,2,1,2,1}=v_\lambda$, and $L'(\lambda)$ is simple.
\epf

\subsection{The family $\If_{25}$}
Recall that $\If_{25}   = \{ \lambda \in \widehat{\Gamma} \mid \lambda_1 = \zeta^{11}, \, \lambda_2=\zeta^{8} \}$.

\begin{lemma}\label{lem:modulo irr caso 25}
If $\lambda \in \If_{25}$, then $\dim L(\lambda)=37$. A basis of $L(\lambda)$ is given by 
\begin{align*}
\mathrm{B}_{25}&= \mathrm{B}'_{25}- \Big( \{ m_{0,3,0,0,e}\mid e \in \I_{0,1}\} \cup \{m_{1,3,c,0,e}, 
m_{1,2,1,0,e} \mid c \in \I_{0,1},e \in \I_{0,2}\}\Big), \text{ where }\\
\mathrm{B}'_{25}&= \{m_{a,b,c,0,e}\mid a,c \in \I_{0,1}, b \in\I_{0,3}, e \in \I_{0,2} \}. \\ 
\end{align*}
\end{lemma}

\pf
Let $w_1= F_1^2 E_{112}E_1^2 v_\lambda$. By Remark \ref{rem:rel E112}, $F_i w_1=0$, $i\in \I_2$. Let $M'=M(\lambda)/\Uc w_1$, so
$\mathrm{B}'_{25}$ is a basis of $M'$. Notice that
$w_2=E_2E_{12}^3v_\lambda$ satisfies $F_1w_2=F_2w_2=0$. 
Hence $\Uc w_2 \twoheadrightarrow L(\mu)$ for $\mu \in \If_{38}$, and there exists $E\in\Uc$ such that
$Ew_2=m_{1,3,1,0,2}$.  Moreover, there exists $F \in \Uc$
such that $Fm_{1,3,1,0,2}=w_2$, and then $\Uc w_2 = \Uc m_{1,3,1,0,2}\simeq L(\mu)$.
Let $L'(\lambda)=M(\lambda)/\Uc w_1 + \Uc w_2$, so $\dim L'(\lambda)=48-11=37$ by Lemma \ref{lem:modulo irr caso 38} and $\mathrm{B}_{25}$ is a basis of $L'(\lambda)$. As in previous cases, $L'(\lambda)$ is simple.
\epf

\subsection{The family $\If_{26}$}
Recall that $\If_{26}   = \{ \lambda \in \widehat{\Gamma} \mid \lambda_1 = \zeta^{5}, \, \lambda_2=\zeta^{8} \}$.

\begin{lemma}\label{lem:modulo irr caso 26}
If $\lambda \in \If_{26}$, then $\dim L(\lambda)=25$. A basis of $L(\lambda)$ is given by
\begin{align*}
\mathrm{B}_{26}&=\{m_{0,b,c,0,e}| c \in \I_{0,1}, b \in\I_{0,3}, e \in \I_{0,2} \} \cup \{m_{1,0,0,0,0}, m_{1,0,0,0,2}\} - \{m_{0,3,1,0,0} \}.
\end{align*}
\end{lemma}

\pf
Let $w_1= F_1^2 E_{112}E_1^2 v_\lambda$, so $F_i w_1=0$, $i\in \I_2$. Let $M'=M(\lambda)/\Uc w_1$. Then 
$\mathrm{B}'_{25}$ as in Lemma \ref{lem:modulo irr caso 26} is a basis of $M'$. Notice that
$w_2=F_2E_2E_{12}v_\lambda$ satisfies $F_1w_2=F_2w_2=0$. 
Hence $\Uc w_2 \twoheadrightarrow L(\mu)$ for $\mu \in \If_{13}$, and there exists $E\in\Uc$ such that
$Ew_2=m_{1,3,1,0,2}$.  Moreover, there exists $F \in \Uc$
such that $Fm_{1,3,1,0,2}=w$, and then $\Uc w_2 = \Uc m_{1,3,1,0,2}\simeq L(\mu)$.
Let $L'(\lambda)=M(\lambda)/\Uc w_1 + \Uc w_2$, so $\dim L'(\lambda)=48-23=25$ by Lemma \ref{lem:modulo irr caso 13}, and $\mathrm{B}_{26}$ is a basis of $L'(\lambda)$. As in previous cases, $L'(\lambda)$ is simple.
\epf

\subsection{The family $\If_{27}$}

Recall that $\If_{27}   = \{ \lambda \in \widehat{\Gamma} \mid \lambda_1 = \zeta^{4}, \, \lambda_2=\zeta^{9} \}$.

\begin{lemma}\label{lem:modulo irr caso 27}
If $\lambda \in \If_{27}$, then $\dim L(\lambda)=35$. A basis of $L(\lambda)$ is given by \begin{align*}
\mathrm{B}_{27}&=\mathrm{B}'_{27}- \{n_{0,0,1,2,2}\}, & &\text{where} & 
\mathrm{B}'_{27}&=\{n_{a,0,c,d,e}\mid a,c \in \I_{0,1},  d,e \in \I_{0,2}\}.
\end{align*}
\end{lemma}

\pf 
Let $w_1=F_2E_{12}E_2v_\lambda$, so $F_i w_1=0$, $i\in \I_2$. Let $M'=M(\lambda)/\Uc w_1$. Then 
$\mathrm{B}'_{27}$ is a basis of $M'$. Notice that
$w_2=E_{11212}E_{112}^2E_1^2v_\lambda$ satisfies $F_1w_2=F_2w_2=0$. 
Hence $\Uc w_2 \twoheadrightarrow L(\mu)$ for $\mu \in \If_{47}$; as also $E_1w_2=E_2w_2=0$, we have that $\Uc w_2 \simeq L(\mu)$.
Let $L'(\lambda)=M(\lambda)/\Uc w_1 + \Uc w_2$, so $\dim L'(\lambda)=36-1=35$ by Lemma \ref{lem:modulo irr caso 47}, and $\mathrm{B}_{27}$ is a basis of $L'(\lambda)$. As in previous cases, $L'(\lambda)$ is simple.
\epf

\subsection{The family $\If_{28}$}

Recall that $\If_{28}   = \{ \lambda \in \widehat{\Gamma} \mid \lambda_1 = \zeta^{9}, \, \lambda_2= \zeta^{4} \}$.

\begin{lemma}\label{lem:modulo irr caso 28}
If $\lambda \in \If_{28}$, then $\dim L(\lambda)=25$. A basis of $L(\lambda)$ is given by
\begin{align*}
\mathrm{B}_{28}&=\mathrm{B}'_{27} - \Big(\{ n_{0,0,1,1,e}, n_{0,0,c,2,e} |
c \in \I_{0,1}, e \in \I_{0,2}\} \cup \{ n_{1,0,1,2,e} \mid e \in \I_{1,2} \} \Big). 
%\text{ where } \\
%\widetilde{\mathrm{B}}_{28}&=\{n_{a,0,c,d,e}| a,c \in \I_{0,1}, d,e \in \I_{0,2} \}
\end{align*}
\end{lemma}

\pf
Let $w_1=F_2E_{12}E_2v_\lambda$, so $F_i w_1=0$, $i\in \I_2$. Let $M'=M(\lambda)/\Uc w_1$. Then 
$\mathrm{B}'_{27}$ is a basis of $M'$. Notice that
$w_2=F_1^2E_1^2E_{112}^2v_\lambda$ satisfies $F_1w_2=F_2w_2=0$. 
Hence $\Uc w_2 \twoheadrightarrow L(\mu)$ for $\mu \in \If_{38}$, and there exists $E\in\Uc$ such that
$Ew_2=m_{1,0,1,2,2}$.  Moreover, there exists $F \in \Uc$
such that $Fm_{1,0,1,2,2}=w_2$, and then $\Uc w_2 = \Uc m_{1,0,1,2,2}\simeq L(\mu)$.
Let $L'(\lambda)=M(\lambda)/\Uc w_1 + \Uc w_2$, so $\dim L'(\lambda)=36-11=25$ by Lemma \ref{lem:modulo irr caso 38}, and $\mathrm{B}_{28}$ is a basis of $L'(\lambda)$. As in previous cases, $L'(\lambda)$ is simple.
\epf

\subsection{The family $\If_{29}$}

Recall that $\If_{29}   = \{ \lambda \in \widehat{\Gamma} \mid \lambda_1 = -1, \, \lambda_2=-1 \}$.

\begin{lemma}\label{lem:modulo irr caso 29}
If $\lambda \in \If_{29}$, then $\dim L(\lambda)=47$. A basis of $L(\lambda)$ is given by
\begin{align*}
\mathrm{B}_{29} &=\mathrm{B}'_{29}- \{m_{1,3,1,0,0}\}, & &\text{where} & 
\mathrm{B}'_{29}&=\{m_{a,b,c,0,e}| a,c \in \I_{0,1}, b \in\I_{0,3}, e \in \I_{0,2} \}.
\end{align*}
\end{lemma}

\pf
Let $w_1=F_1^2E_{112}E_1^2v_\lambda$, so $F_i w_1=0$, $i\in \I_2$. Let $M'=M(\lambda)/\Uc w_1$. Then 
$\mathrm{B}'_{29}$ is a basis of $M'$. Notice that
$w_2=E_2E_{12}^3E_{11212} v_\lambda$ satisfies $F_1w_2=F_2w_2=0$. 
Hence $\Uc w_2 \twoheadrightarrow L(\mu)$ for $\mu \in \If_{47}$; as also $E_1w_2=E_2w_2=0$, we have that $\Uc w_2 \simeq L(\mu)$.
Let $L'(\lambda)=M(\lambda)/\Uc w_1 + \Uc w_2$, so $\dim L'(\lambda)=48-1=47$ by Lemma \ref{lem:modulo irr caso 47}, and $\mathrm{B}_{29}$ is a basis of $L'(\lambda)$. As in previous cases, $L'(\lambda)$ is simple.
\epf

\subsection{The family $\If_{30}$}
Recall that $\If_{30}   = \{ \lambda \in \widehat{\Gamma} \mid \lambda_1 = \zeta^{2}, \, \lambda_2=\zeta^{2} \}$.

\begin{lemma}\label{lem:modulo irr caso 30}
If $\lambda \in \If_{30}$, then $\dim L(\lambda)=37$. A basis of $L(\lambda)$ is given by
\begin{align*}
\mathrm{B}_{30}&=\mathrm{B}'_{29} 
- \{ m_{1,b,c,0,e}\mid  c \in \I_{0,1}, b \in\I_{2,3}, e \in \I_{0,2}, (b,c,e) \neq (3,1,2) \}. 
%\text{ where }
%\\ \widetilde{\mathrm{B}}_{30}&=\{m_{a,b,c,0,e}| a,c \in \I_{0,1}, b \in\I_{0,3}, e \in \I_{0,2} \}. 
\end{align*}
\end{lemma}
\pf
Let $w_1=F_1^2E_{112}E_1^2v_\lambda$, so $F_i w_1=0$, $i\in \I_2$. Let $M'=M(\lambda)/\Uc w_1$. Then 
$\mathrm{B}'_{29}$ is a basis of $M'$. Notice that
$w_2=E_2E_{12}^2v_\lambda$ satisfies $F_1w_2=F_2w_2=0$. 
Hence $\Uc w_2 \twoheadrightarrow L(\mu)$ for $\mu \in \If_{38}$, and there exists $E\in\Uc$ such that
$Ew_2=m_{1,3,1,0,2}$.  Moreover, there exists $F \in \Uc$
such that $Fm_{1,3,1,0,2}=w_2$, and then $\Uc w_2 = \Uc m_{1,3,1,0,2}\simeq L(\mu)$.
Let $L'(\lambda)=M(\lambda)/\Uc w_1 + \Uc w_2$, so $\dim L'(\lambda)=48-11=37$ by Lemma \ref{lem:modulo irr caso 38}, and $\mathrm{B}_{30}$ is a basis of $L'(\lambda)$. As in previous cases, $L'(\lambda)$ is simple.
\epf

\subsection{The family $\If_{31}$}
Recall that $\If_{31}   = \{ \lambda \in \widehat{\Gamma} \mid \lambda_1 = -1, \, \lambda_2= \zeta^{10} \}$.

\begin{lemma}\label{lem:modulo irr caso 31}
If $\lambda \in \If_{31}$, then $\dim L(\lambda)=61$. A basis of $L(\lambda)$ is given by
\begin{align*}
\mathrm{B}_{31} &= \mathrm{B}'_{31} - \big( \{ n_{0,0,0,2,e}\mid e\in \I_{0,1} \}
\cup \{n_{0,0,1,1,e},n_{0,0,1,2,e}, n_{0,1,1,2,e}|  e \in \I_{0,2}\} \big), \quad \text{where} \\
\mathrm{B}'_{31}&=\{n_{a,b,c,d,e}| a,b,c \in \I_{0,1}, d,e \in \I_{0,2} \}.
\end{align*}
\end{lemma}

\pf
Let $w_1= F_2 E_2E_{12}^2v_\lambda$. By Remark \ref{rem:rel E12 cuadrado}, $F_i w_1=0$, $i\in \I_2$. Let $M'=M(\lambda)/\Uc w_1$, so
$\mathrm{B}'_{31}$ is a basis of $M'$. Notice that
\begin{align*}
w_2&=n_{0,0,0,2,1}+\frac{q_{21}}{3}\zeta(1+\zeta^3)(1+\zeta^2) \big(n_{0,0,1,0,2}+\zeta^4n_{0,1,0,1,2} \big)
\end{align*}
satisfies $F_1w_2=F_2w_2=0$. 
Hence $\Uc w_2 \twoheadrightarrow L(\mu)$ for $\mu \in \If_{18}$, and there exists $E\in\Uc$ such that
$Ew_2=n_{1,1,1,2,2}$.  Moreover, there exists $F \in \Uc$
such that $Fn_{1,1,1,2,2}=w_2$, and then $\Uc w_2 = \Uc n_{1,1,1,2,2}\simeq L(\mu)$.
Let $L'(\lambda)=M(\lambda)/\Uc w_1 + \Uc w_2$, so $\dim L'(\lambda)=72-11=61$ by Lemma \ref{lem:modulo irr caso 18}, and $\mathrm{B}_{31}$ is a basis of $L'(\lambda)$. As in previous cases, $L'(\lambda)$ is simple.
\epf

\subsection{The family $\If_{32}$}
Recall that $\If_{32}   = \{ \lambda \in \widehat{\Gamma} \mid \lambda_1 =\zeta^{10}, \, \lambda_2=-1 \}$.

\begin{lemma}\label{lem:modulo irr caso 32}
If $\lambda \in \If_{32}$, then $\dim L(\lambda)=61$. A basis of $L(\lambda)$ is given by \begin{align*}
\mathrm{B}_{32} &=  \mathrm{B}'_{31} 
- \Big( \{ n_{a,b,1,d,2}\mid a, b \in \I_{0,1},d \in \I_{1,2} \} \cup \{ n_{0,0,1,0,2}, n_{1,0,1,0,2}, n_{1,0,0,2,2}\} \Big).
\end{align*}
\end{lemma}

\pf
Let $w_1= F_2 E_2E_{12}^2v_\lambda$. By Remark \ref{rem:rel E12 cuadrado}, $F_i w_1=0$, $i\in \I_2$. Let $M'=M(\lambda)/\Uc w_1$, so
$\mathrm{B}'_{31}$ is a basis of $M'$. 
Moreover $u=n_{1,1,1,2,2}\in V_{10\alpha_1+6\alpha_2}$ satisfies that $E_1u=E_2u=0$, $g_1\gpl_1 u=u$,
$g_2\gpl_2 u=\zeta^8 u$, so $(\Uc w)^\varphi\twoheadrightarrow L(\nu)$, $\nu\in\If_{12}$. Also $\Uc u$ is a proper submodule. Set $L'(\lambda)=M(\lambda)/\Uc w_1+\Uc u$.  By Lemma \ref{lem:modulo irr caso 12},
$$ 61= \dim L(\lambda)\leq \dim L'(\lambda)=\dim W-\dim \Uc w\leq \dim W-\dim L(\nu)=61, $$
so $L(\lambda)= L'(\lambda)$ and $\Uc w \simeq L(\nu)^\varphi$. In particular $w_2:=F_2F_{11212}F_{112}u\neq 0$, $F_iw_2=0$ and $\Uc w_2=\Uc u$.
Moreover $\mathrm{B}_{32}$ is a basis of $L(\lambda)$.
\epf

\subsection{The family $\If_{33}$}
Recall that $\If_{33}   = \{ \lambda \in \widehat{\Gamma} \mid \lambda_1 = \zeta^{2}, \, \lambda_2= -1 \}$.

\begin{lemma}\label{lem:modulo irr caso 33}
If $\lambda \in \If_{33}$, then $\dim L(\lambda)=71$. A basis of $L(\lambda)$ is given by \begin{align*}
\mathrm{B}_{33}&=  \{m_{a,b,c,d,e} \mid a,c, d \in \I_{0,1},  b,e \in \I_{0,2}
\} \cup \{ m_{1,3,0,0,0}\} - \{ m_{0,0,1,0,0},m_{1,2,0,1,2}\}.
\end{align*}
\end{lemma}

\pf
Let $w_1=F_1^2E_{112}^2E_1^2v_\lambda$. By Remark \ref{rem:rel E112 cuadrado}, $F_1w_1=F_2w_1=0$. By a direct computation, $\Uc w_1 \simeq L(\mu)$, with $\mu\in \If_{23}$, and $\mathrm{B}' =\{ m_{a,b,c,d,e}\mid d \neq 2\} \cup \{ m_{0,0,0,2,2}\}$ is a basis of $W'= M(\lambda)/\Uc w_1$.
%Indeed $ m_{a,b,c,2,e}$ appears with non-zero coefficient in $E_1^eE_{11212}E_{12}^bE_2^aw_1$ if $(a,b,c,e) \neq (0,0,0,2)$, but $E_1^2 w_1=0$ by direct computation, so $\mathrm{B}'$ is linearly independent in $W'$. It is a basis since
%$\dim W'= 144- \dim \Uc w_1 \geq 144 - \dim L(\mu) = 97$.
Now $\Uc m_{0,0,0,2,2} = \ku m_{0,0,0,2,2}$ in $W'$, so $\mathrm{B}= \{m_{a,b,c,d,e}\mid d\neq 2\}$ is a basis of $M'= W'/\ku m_{0,0,0,2,2}$.

Let $w_2=F_1^2F_{112}^2E_{11212}E_{112}^2E_1^2 v_\lambda$. By Remark \ref{rem:l1 a la 3 l2 a la 2}, $F_iw_2=0$, $i\in \I_2$, and $\Uc w_2\twoheadrightarrow L(\mu)$, with $\mu\in \If_{14}$, and there exists $E\in\Uc$ such that
$Ew_2=m_{1,3,1,1,2}$.  Moreover, there exists $F \in \Uc$
such that $Fm_{1,3,1,1,2}=w_2$, and then $\Uc w_2 = \Uc m_{1,3,1,1,2}\simeq L(\mu)$.
Let $L'(\lambda)=M(\lambda)/\Uc w_1 +\Uc m_{0,0,0,2,2} + \Uc w_2$, so $\dim L'(\lambda)=96-25=71$ by Lemma \ref{lem:modulo irr caso 14}, and $\mathrm{B}_{33}$ is a basis of $L'(\lambda)$. As in previous cases, $L'(\lambda)$ is simple.
\epf

\subsection{The family $\If_{34}$}
Recall that $\If_{34}   = \{ \lambda \in \widehat{\Gamma} \mid \lambda_1 = \zeta^{4}, \, \lambda_2= \zeta^{3} \}$.

\begin{lemma}\label{lem:modulo irr caso 34}
If $\lambda \in \If_{34}$, then $\dim L(\lambda)=71$. A basis of $L(\lambda)$ is given by
\begin{align*}
\mathrm{B}_{34}= & \{n_{a,b,c,d,e}| a,c, d \in \I_{0,1},  b,e \in \I_{0,2}\} \cup
\{n_{0,0,0,2,e} \mid e \in \I_{0,2} \} \\ &- \big(\{n_{0,0,1,0,e}| e \in \I_{0,2}\} \cup \{n_{0,1,1,1,0} \}\big).
\end{align*}
\end{lemma}

\pf
Let $w_1=F_2E_{12}^3E_2v_\lambda$. By Remark \ref{rem:rel E12 cubo}, $F_1w_1=F_2w_1=0$. By a direct computation, $\Uc w_1 \simeq L(\mu)$, with $\mu\in \If_{36}$, and $\mathrm{B}' = \mathrm{B}'_{35} \cup \{n_{1,3,0,0,0}\}$ is a basis of $W'= M(\lambda)/\Uc w_1$.
Now $\Uc n_{1,3,0,0,0} = \ku n_{1,3,0,0,0}$ in $W'$, so $\mathrm{B}'_{35}$ is a basis of $M'= W'/\ku n_{1,3,0,0,0}$.

Let $w_2=F_1^2F_{112}^2E_{11212}E_{112}^2E_1^2 v_\lambda$. By Remark \ref{rem:l1 a la 3 l2 a la 2}, $F_iw_2=0$, $i\in \I_2$, and $\Uc w_2\twoheadrightarrow L(\mu)$, with $\mu\in \If_{37}$, and there exists $E\in\Uc$ such that
$Ew_2=n_{1,2,1,2,2}$.  Moreover, there exists $F \in \Uc$
such that $Fn_{1,2,1,2,2}=w_2$, and then $\Uc w_2 = \Uc n_{1,2,1,2,2}\simeq L(\mu)$.
Let $L'(\lambda)=M(\lambda)/\Uc w_1 +\Uc n_{1,2,1,2,2} + \Uc w_2$, so $\dim L'(\lambda)=108-37=71$ by Lemma \ref{lem:modulo irr caso 37}, and $\mathrm{B}_{34}$ is a basis of $L'(\lambda)$. As in previous cases, $L'(\lambda)$ is simple.
\epf

\subsection{The family $\If_{35}$}
Recall that $\If_{35}   = \{ \lambda \in \widehat{\Gamma} \mid \lambda_1 = \zeta^{3}, \, \lambda_2=\zeta^{4} \}$.

\begin{lemma}\label{lem:modulo irr caso 35}
If $\lambda \in \If_{35}$, then $\dim L(\lambda)=85$. A basis of $L(\lambda)$ is given by
\begin{align*}
\mathrm{B}_{35}= & \mathrm{B}'_{35} 
- \Big(\{ n_{0,b,c,2,e} | c \in \I_{0,1},  b, e \in \I_{0,2}\} \cup \{n_{1,2,1,2,2}, n_{1,0,0,2,2}, n_{1,0,1,2,e}| e \in \I_{0,2}\}\Big) \text{ where }
\\
\mathrm{B}'_{35}= & \{n_{a,b,c,d,e} \mid a,c \in \I_{0,1},  b, d,e \in \I_{0,2}\}
\end{align*}
\end{lemma}

\pf
Let $w_1=F_2E_2E_{12}^3v_\lambda$, so $F_i w_1=0$, $i\in \I_2$. Let $M'=M(\lambda)/\Uc w_1$. Then 
$\mathrm{B}'_{35}$ is a basis of $M'$. Notice that
$w_2=F_1^2E_{112}^2E_1^2v_\lambda$ satisfies $F_1w_2=F_2w_2=0$. 
Hence $\Uc w_2 \twoheadrightarrow L(\mu)$ for $\mu \in \If_{44}$, and there exists $E\in\Uc$ such that
$Ew_2=n_{1,2,1,2,2}$.  Moreover, there exists $F \in \Uc$
such that $Fn_{1,2,1,2,2}=w_2$, and then $\Uc w_2 = \Uc n_{1,2,1,2,2}\simeq L(\mu)$.
Let $L'(\lambda)=M(\lambda)/\Uc w_1 + \Uc w_2$, so $\dim L'(\lambda)=108-23=85$ by Lemma \ref{lem:modulo irr caso 44}, and $\mathrm{B}_{35}$ is a basis of $L'(\lambda)$. As in previous cases, $L'(\lambda)$ is simple.
\epf

\subsection{The family $\If_{36}$}
Recall that $\If_{36}   = \{ \lambda \in \widehat{\Gamma} \mid \lambda_1 = \zeta, \, \lambda_2= 1 \}$.

\begin{lemma}\label{lem:modulo irr caso 36}
If $\lambda \in \If_{36}$, then $\dim L(\lambda)=35$. A basis of $L(\lambda)$ is given by $\mathrm{B}_{36}=$
\begin{align*}
 &\{n_{0,b,0,d,e},n_{0,0,1,2,e},n_{0,0,1,0,e}|  b \in\I_{0,3}, d,e \in \I_{0,2}\} - \{n_{0,1,0,1,e}, n_{0,2,0,2,e},n_{0,1,0,0,2} | e \in \I_{0,2}\}.
\end{align*}
\end{lemma}

\pf
Let $w_1=\widetilde{n}_{1,0,0,0,0}$, $w_2 =E_1^2E_{12}v_\lambda$. Then $W(\lambda)=\Uc w_1$ and $F_1w_2=F_2w_2=0$.
Set $M'=M(\lambda)/W_2(\lambda)$, $u=\widetilde{n}_{0,3,1,2,2}$.
Hence $\Uc w_2\twoheadrightarrow L(\mu)$ for $\mu \in \If_{15}$, and there exists $E\in\Uc$ such that
$Ew_2=u$.  Moreover, there exists $F \in \Uc$
such that $Fu=w_2$, and then $\Uc w_2 = \Uc u\simeq L(\mu)$.
Let $L'(\lambda)=M(\lambda)/\Uc w_2 + W(\lambda)$, so $\dim L'(\lambda)=72-37=35$ by Lemma \ref{lem:modulo irr caso 15}, and $\mathrm{B}_{36}$ is a basis of $L'(\lambda)$. As in previous cases, $L'(\lambda)$ is simple.
\epf

\subsection{The family $\If_{37}$}
Recall that $\If_{37}   = \{ \lambda \in \widehat{\Gamma} \mid \lambda_1 = \zeta^{2}, \, \lambda_2= 1 \}$.

\begin{lemma}\label{lem:modulo irr caso 37}
If $\lambda \in \If_{37}$, then $\dim L(\lambda)=37$. A basis of $L(\lambda)$ is given by
\begin{align*}
\mathrm{B}_{37} &= \{n_{0,b,0,d,e},n_{0,0,1,0,0}, n_{0,3,1,0,e}|  b \in\I_{0,3}, d,e \in \I_{0,2}\}-\{ n_{0,3,0,2,e} | e \in \I_{0,2} \}.
\end{align*}
\end{lemma}

\pf
Let $w_1=\widetilde{n}_{1,0,0,0,0}$, $w_2 =\widetilde{n}_{0,1,0,1,1}-\zeta \widetilde{n}_{0,2,0,0,2}- \zeta^{10}(1-\zeta)^2 \widetilde{n}_{0,0,1,0,1}$. Then $W(\lambda)=\Uc w_1$ and $F_1w_2=F_2w_2=0$.
Set $M'=M(\lambda)/W_2(\lambda)$, $u=\widetilde{n}_{0,3,1,2,2}$.
Hence $\Uc w_2\twoheadrightarrow L(\mu)$ for $\mu \in \If_{19}$, and there exists $E\in\Uc$ such that
$Ew_2=u$.  Moreover, there exists $F \in \Uc$
such that $Fu=w_2$, and then $\Uc w_2 = \Uc u\simeq L(\mu)$.
Let $L'(\lambda)=M(\lambda)/\Uc w_2 + W(\lambda)$, so $\dim L'(\lambda)=72-35=37$ by Lemma \ref{lem:modulo irr caso 19}, and $\mathrm{B}_{37}$ is a basis of $L'(\lambda)$. As in previous cases, $L'(\lambda)$ is simple.
\epf

\subsection{The family $\If_{38}$}
Recall that $\If_{38}   = \{ \lambda \in \widehat{\Gamma} \mid \lambda_1 = \zeta^{3}, \, \lambda_2= 1 \}$.

\begin{lemma}\label{lem:modulo irr caso 38}
If $\lambda \in \If_{38}$, then $\dim L(\lambda)=11$. A basis of $L(\lambda)$ is given by
\begin{align*}
\mathrm{B}_{38}= &\{n_{0,b,c,0,e}| b,c \in \I_{0,1}, e \in \I_{0,2}\}-\{n_{0,1,1,0,2}\}.
\end{align*}
The action of $E_i$, $F_i$, $i\in \I_2$ is described in Table \ref{tab:caso 38}.
\end{lemma}

\begin{table}[h]
	\caption{Simple modules for $\lambda \in \If_{38}$}\label{tab:caso 38}
	\begin{center}
		\begin{tabular}{| c | c | c | c | c |} \hline
			$w$ & $E_1\cdot w$ & $E_2\cdot w$ & $\lambda(g_1^{-1})F_1\cdot w$ & $\lambda(g_2^{-1})F_2\cdot w$ \\ \hline
			$v_{0,0}$ & $v_{1,0}$ & $0$ & $0$ & $0$ \\ \hline
			$v_{1,0}$ & $v_{2,0}$ & $\zeta^7q_{21}v_{1,1}$ & $(1-\zeta^3)v_{0,0}$ & $0$ \\ \hline
			$v_{2,0}$ & $0$ & $\zeta^8q_{21}^2(1+\zeta^3)v_{2,1}$ & $\zeta^7(1+\zeta)v_{1,0}$ & $0$\\ \hline
			$v_{1,1}$ & $v_{2,1}$ & $0$ & $0$ & $(\zeta^{11}-1)v_{1,0}$ \\ \hline
			$v_{2,1}$ & $v_{3,1}$ & $0$ & $q_{12}\zeta^8v_{1,1}$ & $(\zeta^{11}-1)v_{2,0}$\\  \hline
			$v_{3,1}$ & $0$ & $q_{21}^2\zeta v_{3,2}$ & $q_{12}\zeta^2v_{2,1}$ & $0$\\ \hline
			$v_{3,2}$ & $v_{4,3}$ & $0$ & $0$ & $q_{21}\zeta^{11}(1-\zeta^3)v_{3,1}$\\ \hline
			$v_{4,2}$ & $v_{5,2}$ & $q_{21}^2\zeta^{10}v_{4,3}$ & $q_{12}^2(\zeta^{11}-1)v_{3,2}$  & $0$ \\  \hline
			$v_{5,2}$ & $0$ & $q_{21}^3(3)_{\zeta}v_{5,3}$ & $q_{12}^2\zeta^8(1+\zeta)v_{4,2}$  & $0$ \\ \hline
			$v_{4,3}$ & $v_{5,3}$ & $0$ & $0$ & $q_{21}^2\zeta^{10}(3)_{\zeta^{11}}v_{4,2}$\\ \hline
			$v_{5,3}$ & $0$ & $0$ & $q_{12}^3\zeta^8(1+\zeta^2)v_{4,3}$ & $q_{21}^2\zeta^{10}(3)_{\zeta^{11}}v_{5,2}$\\ \hline
		\end{tabular}
	\end{center}
\end{table}

\pf
Let $w_1=\widetilde{n}_{1,0,0,0,0}$, $w_2 =F_1^2E_{112}E_1^2v_\lambda$. Then $W(\lambda)=\Uc w_1$ and $F_1w_2=F_2w_2=0$.
Let $L'(\lambda)=M(\lambda)/\Uc w_2 +W(\lambda)$. We label the elements of $\mathrm{B}_{38}$ as follows:
\begin{align*}
v_{0,0}&=n_{0,0,0,0,0}, & v_{1,1}&=n_{0,1,0,0,0}, &
v_{3,2}&=n_{0,0,1,0,0}, & v_{4,3}&=n_{0,1,1,0,0}, \\
v_{1,0}&=n_{0,0,0,0,1}, & v_{2,1}&=n_{0,1,0,0,1}, &
v_{4,2}&=n_{0,0,1,0,1}, & v_{5,3}&=n_{0,1,1,0,1}, \\
v_{2,0}&=n_{0,0,0,0,2}, & v_{3,1}&=n_{0,1,0,0,2}, &
v_{5,2}&=n_{0,0,1,0,2}.
\end{align*}
We check that the action of $E_k$, $F_k$ on $v_{i,j}$ is given by Table \ref{tab:caso 38} and $L'(\lambda)$ is spanned by $\mathrm{B}_{38}$.
Moreover there exists $F\in\Uc^-$ such that $Fv_{5,3}=v_\lambda$, and for each pair $(i,j)$ there is $E_{i,j}\in\Uc_{(5-i)\alpha_1+(3-j)\alpha_2}$ such that $E_{i,j}v_{i,j}=v_{5,3}$. Thus $L'(\lambda)$ is simple.
\epf

\subsection{The family $\If_{39}$}
Recall that $\If_{39}   = \{ \lambda \in \widehat{\Gamma} \mid \lambda_1 = \zeta^{4}, \, \lambda_2= 1 \}$.

\begin{lemma}\label{lem:modulo irr caso 39}
If $\lambda \in \If_{39}$, then $\dim L(\lambda)=61$. A basis of $L(\lambda)$ is given by
\begin{align*}
\mathrm{B}_{39}= &\{ n_{0,b,c,d,e} | c \in \I_{0,1}, b \in\I_{0,3}, d,e \in \I_{0,2}\} 
\\& -\Big( \{n_{0,3,c,2,e}, n_{0,2,1,2,e} | c \in \I_{0,1},e \in \I_{0,2}\} \cup \{n_{0,2,0,2,e}\mid e\in \I_{1,2}\}\Big).
\end{align*}
\end{lemma}
\pf
Let $w_1=\widetilde{n}_{1,0,0,0,0}$, $u=\widetilde{n}_{0,3,1,2,2}$, $w_2 =F_1F_{11212}F_{12}^2 u$. Then 
$W_2(\lambda)=\Uc w_1$.

Let $M'=M(\lambda)/W(\lambda)$, so $E_1u=E_2u=0$ in $M'$, and $(\mathcal U u)^{\varphi}\twoheadrightarrow L(\nu)$ for some  $\nu\in \If_{38}$; thus $w_2\neq 0$.
By direct computation, $F_iw_2=0$, $i\in\I_2$, so $\mathcal{U} w_2$ projects over a simple module $L(\mu)$, for $\mu \in \If_{18}$. From here, $\mathcal{U} w_2\simeq L(\mu)$.

Let $L'(\lambda)=M(\lambda)/W_1(\lambda)+\Uc w_2$. Then $\dim L'(\lambda)=61$ by Lemma \ref{lem:modulo irr caso 18}, and $\mathrm{B}_{39}$ is a basis of $L'(\lambda)$. There exists $F$ such that $F u=v_\lambda$, and $L'(\lambda)$ is simple.
\epf

\subsection{The family $\If_{40}$}
Recall that $\If_{40}   = \{ \lambda \in \widehat{\Gamma} \mid \lambda_1 = \zeta^{5}, \, \lambda_2= 1 \}$.

\begin{lemma}\label{lem:modulo irr caso 40}
If $\lambda \in \If_{40}$, then $\dim L(\lambda)=35$. A basis of $L(\lambda)$ is given by
\begin{align*}
\mathrm{B}_{40}= &\{ n_{0,b,c,0,e} | c \in \I_{0,1}, b \in\I_{0,3}, e \in \I_{0,2}\} 
\cup \{n_{0,b,c,1,e} \mid c \in \I_{0,1}, b \in\I_{0,1}, e \in \I_{0,2} \}
\\ & \cup \{ n_{0,3,0,2,e}\mid e \in \I_{0,1} \} -\{n_{0,3,1,0,e} | e \in \I_{0,2} \} .
 \end{align*}
\end{lemma}
\pf
Let $w_1=\widetilde{n}_{1,0,0,0,0}$, $w_2 =F_1^2E_{112}^2E_1^2v_\lambda$. Then $W(\lambda)=\Uc w_1$ and $F_1w_2=F_2w_2=0$.
Set $M'=M(\lambda)/W_2(\lambda)$, $u=\widetilde{n}_{0,3,1,2,2}$.
Hence $\Uc w_2\twoheadrightarrow L(\mu)$ for $\mu \in \If_{25}$, and there exists $E\in\Uc$ such that
$Ew_2=u$.  Moreover, there exists $F \in \Uc$
such that $Fu=w_2$, and then $\Uc w_2 = \Uc u\simeq L(\mu)$.
Let $L'(\lambda)=M(\lambda)/\Uc w_2 + W(\lambda)$, so $\dim L'(\lambda)=72-37=35$ by Lemma \ref{lem:modulo irr caso 25}, and $\mathrm{B}_{40}$ is a basis of $L'(\lambda)$. As in previous cases, $L'(\lambda)$ is simple.
\epf

\subsection{The family $\If_{41}$}
Recall that $\If_{41}   = \{ \lambda \in \widehat{\Gamma} \mid \lambda_1 = -1, \, \lambda_2= 1 \}$.

\begin{lemma}\label{lem:modulo irr caso 41}
If $\lambda \in \If_{41}$, then $\dim L(\lambda)=37$. A basis of $L(\lambda)$ is given by
\begin{align*}
\mathrm{B}_{41} &= \{n_{0,b,c,d,0}| c \in \I_{0,1}, b, d \in \I_{0,2} \}  
\cup  \{n_{0,b,c,d,e}| c, b \in \I_{0,1}, d \in \I_{0,2},  e\in \I_{1,2}\} \\
& - \{n_{0,1,c,d,2},n_{0,0,1,2,2}| c \in \I_{0,1} d\in \I_{1,2} \}.
\end{align*}
\end{lemma}
\pf
Let $w_1=\widetilde{n}_{1,0,0,0,0}$, $w_2 =F_1^2F_{112}^2E_{11212}E_{112}^2E_1^2v_\lambda$. Then $W(\lambda)=\Uc w_1$ and $F_1w_2=F_2w_2=0$.
Set $M'=M(\lambda)/W_2(\lambda)$, $u=\widetilde{n}_{0,3,1,2,2}$.
Hence $\Uc w_2\twoheadrightarrow L(\mu)$ for $\mu \in \If_{27}$, and there exists $E\in\Uc$ such that
$Ew_2=u$.  Moreover, there exists $F \in \Uc$
such that $Fu=w_2$, and then $\Uc w_2 = \Uc u\simeq L(\mu)$.
Let $L'(\lambda)=M(\lambda)/\Uc w_2 + W(\lambda)$, so $\dim L'(\lambda)=72-35=37$ by Lemma \ref{lem:modulo irr caso 27}, and $\mathrm{B}_{41}$ is a basis of $L'(\lambda)$. As in previous cases, $L'(\lambda)$ is simple.
\epf

\subsection{The family $\If_{42}$}
Recall that $\If_{42}   = \{ \lambda \in \widehat{\Gamma} \mid \lambda_1 = \zeta^{7}, \, \lambda_2= 1 \}$.

\begin{lemma}\label{lem:modulo irr caso 42}
If $\lambda \in \If_{42}$, then $\dim L(\lambda)=71$. A basis of $L(\lambda)$ is given by
\begin{align*}
\mathrm{B}_{42}= &\{n_{0,b,c,d,e}| c \in \I_{0,1}, b \in\I_{0,3}, d,e \in \I_{0,2}, (b,c,d,e) \neq (3,1,2,2) \}. 
\end{align*}
\end{lemma}

\pf
Let $w_1=\widetilde{n}_{1,0,0,0,0}$, $w_2 =\widetilde{n}_{0,3,1,2,2}$. Then $W(\lambda)=\Uc w_1$ and $F_1w_2=F_2w_2=E_1w_2=E_2w_2=0$, so $\Uc w_2 \simeq L(\mu)$ for $\mu\in\If_{47}$.
Let $L'(\lambda)=M(\lambda)/W(\lambda) + \Uc w_2$, so $\mathrm{B}_{42}$ is a basis of $L'(\lambda)$. There exists $F\in\Uc^-$ such that
$F n_{0,3,1,2,1}=v_\lambda$.
If $n_{0,b,c,d,e} \in \mathrm{B}_{42}$, then $E_1^{1-e}E_{112}^{2-d}E_{11212}^{1-c}E_{12}^{3-b} n_{0,b,c,d,e} \in \kut n_{0,3,1,2,1}$, so $L'(\lambda)$ is simple.
\epf

\subsection{The family $\If_{43}$}
Recall that $\If_{43}   = \{ \lambda \in \widehat{\Gamma} \mid \lambda_1 = \zeta^{8}, \, \lambda_2= 1 \}$.

\begin{lemma}\label{lem:modulo irr caso 43}
If $\lambda \in \If_{43}$, then $\dim L(\lambda)=25$. A basis of $L(\lambda)$ is given by
\begin{align*}
\mathrm{B}_{43} &= \{n_{0,b,c,d,e}| c,e \in \I_{0,1}, b \in\I_{0,3}, d \in \I_{0,2}\} -\Big( \{n_{0,2,1,2,0}\} \\
& \qquad \cup \{n_{0,b,c,d,1}|c \in \I_{0,1}, b \in\I_{1,3}, d \in \I_{0,2}\} \cup \{n_{0,3,c,d,0}|
c \in \I_{0,1}, d \in \I_{1,2}\} \Big).
\end{align*}
\end{lemma}

\pf
Let $w_1=\widetilde{n}_{1,0,0,0,0}$, $w_2 =E_1^2v_\lambda$. Then $W(\lambda)=\Uc w_1$ and $F_1w_2=F_2w_2=0$.
Set $M'=M(\lambda)/W_2(\lambda)$, $u=\widetilde{n}_{0,3,1,2,2}$.
Hence $\Uc w_2\twoheadrightarrow L(\mu)$ for $\mu \in \If_{17}$, and there exists $E\in\Uc$ such that
$Ew_2=u$.  Moreover, there exists $F \in \Uc$
such that $Fu=w_2$, and then $\Uc w_2 = \Uc u\simeq L(\mu)$.
Let $L'(\lambda)=M(\lambda)/\Uc w_2 + W(\lambda)$, so $\dim L'(\lambda)=72-47=25$ by Lemma \ref{lem:modulo irr caso 17}, and $\mathrm{B}_{43}$ is a basis of $L'(\lambda)$. As in previous cases, $L'(\lambda)$ is simple.
\epf

\subsection{The family $\If_{44}$}
Recall that $\If_{44}   = \{ \lambda \in \widehat{\Gamma} \mid \lambda_1 = \zeta^{9}, \, \lambda_2= 1 \}$.

\begin{lemma}\label{lem:modulo irr caso 44}
If $\lambda \in \If_{44}$, then $\dim L(\lambda)=23$. A basis of $L(\lambda)$ is given by
\begin{align*}
\mathrm{B}_{44}= &\{n_{0,b,0,d,e}|  b \in\I_{0,3}, d \in \I_{0,2}, e \in \I_{0,1}\}
\cup \{n_{0,0,0,0,2}\} - \{ n_{0,3,0,1,1},n_{0,3,0,2,1}\}. 
\end{align*}
\end{lemma}
\pf
Let $w_1=\widetilde{n}_{1,0,0,0,0}$, $w_2 =\zeta^4 \widetilde{n}_{0,0,0,1,1} +\widetilde{n}_{0,1,0,0,2}$. Then $W(\lambda)=\Uc w_1$ and $F_1w_2=F_2w_2=0$.
Set $M'=M(\lambda)/W_2(\lambda)$, $u=\widetilde{n}_{0,3,1,2,2}$.
Hence $\Uc w_2\twoheadrightarrow L(\mu)$ for $\mu \in \If_{22}$, and there exists $E\in\Uc$ such that
$Ew_2=u$.  Moreover, there exists $F \in \Uc$
such that $Fu=w_2$, and then $\Uc w_2 = \Uc u\simeq L(\mu)$.
Let $L'(\lambda)=M(\lambda)/\Uc w_2 + W(\lambda)$, so $\dim L'(\lambda)=72-49=23$  by Lemma \ref{lem:modulo irr caso 22}, and $\mathrm{B}_{44}$ is a basis of $L'(\lambda)$. As in previous cases, $L'(\lambda)$ is simple.
\epf

\subsection{The family $\If_{45}$}
Recall that $\If_{45}   = \{ \lambda \in \widehat{\Gamma} \mid \lambda_1 = \zeta^{10}, \, \lambda_2= 1 \}$.

\begin{lemma}\label{lem:modulo irr caso 45}
If $\lambda \in \If_{45}$, then $\dim L(\lambda)=49$. A basis of $L(\lambda)$ is given by
\begin{align*}
\mathrm{B}_{45}= & \{n_{0,b,c,d,e} | c \in \I_{0,1}, b \in\I_{0,3}, d,e \in \I_{0,2} \}  
\\ &- \big(\{n_{0,b,c,2,e}| c \in \I_{0,1}, b \in\I_{1,3},e \in \I_{0,2} \} \cup \{n_{0,0,1,2,e} | e \in \I_{0,2}\} \cup \{n_{0,0,1,0,2}, n_{0,3,1,1,2}\}\big).
\end{align*}
\end{lemma}

\pf
Let $w_1=\widetilde{n}_{1,0,0,0,0}$, $w_2 =n_{0,1,0,1,2} - \zeta^{11}(3)_{\zeta^7}n_{0,0,1,0,2}$. Then $W(\lambda)=\Uc w_1$ and $F_1w_2=F_2w_2=0$.
Set $M'=M(\lambda)/W_2(\lambda)$, $u=\widetilde{n}_{0,3,1,2,2}$.
Hence $\Uc w_2\twoheadrightarrow L(\mu)$ for $\mu \in \If_{13}$, and there exists $E\in\Uc$ such that
$Ew_2=u$.  Moreover, there exists $F \in \Uc$
such that $Fu=w_2$, and then $\Uc w_2 = \Uc u\simeq L(\mu)$.
Let $L'(\lambda)=M(\lambda)/\Uc w_2 + W(\lambda)$, so $\dim L'(\lambda)=72-23=49$ by Lemma \ref{lem:modulo irr caso 13}, and $\mathrm{B}_{45}$ is a basis of $L'(\lambda)$. As in previous cases, $L'(\lambda)$ is simple.
\epf

\subsection{The family $\If_{46}$}
Recall that $\If_{46}   = \{ \lambda \in \widehat{\Gamma} \mid \lambda_1 = \zeta^{11}, \, \lambda_2= 1 \}$.

\begin{lemma}\label{lem:modulo irr caso 46}
If $\lambda \in \If_{46}$, then $\dim L(\lambda)=47$. A basis of $L(\lambda)$ is
\begin{align*}
\mathrm{B}_{46}&=\{n_{0,b,c,d,e},|c,d \in \I_{0,1}, b \in\I_{0,3}, e \in \I_{0,2}\} \cup \{n_{0,1,0,2,0},n_{0,3,1,2,0}\}\\ &-\{n_{0,1,1,0,2},n_{0,3,0,0,1},n_{0,1,1,0,1}\}.
\end{align*}
\end{lemma}
\pf
Let $w_1=\widetilde{n}_{1,0,0,0,0}$, $w_2 =F_1^2E_{112}^2E_1^2v_\lambda$. Then $W(\lambda)=\Uc w_1$ and $F_1w_2=F_2w_2=0$.
Set $M'=M(\lambda)/W_2(\lambda)$, $u=\widetilde{n}_{0,3,1,2,2}$.
Hence $\Uc w_2\twoheadrightarrow L(\mu)$ for $\mu \in \If_{26}$, and there exists $E\in\Uc$ such that
$Ew_2=u$.  Moreover, there exists $F \in \Uc$
such that $Fu=w_2$, and then $\Uc w_2 = \Uc u\simeq L(\mu)$.
Let $L'(\lambda)=M(\lambda)/\Uc w_2 + W(\lambda)$, so $\dim L'(\lambda)=72-25=47$ by Lemma \ref{lem:modulo irr caso 26}, and $\mathrm{B}_{46}$ is a basis of $L'(\lambda)$. As in previous cases, $L'(\lambda)$ is simple.
\epf

\subsection{The family $\If_{47}$}
Recall that $\If_{47}   = \{ \lambda \in \widehat{\Gamma} \mid \lambda_1 = 1, \, \lambda_2= 1 \}$.

\begin{lemma}\label{lem:modulo irr caso 47}
If $\lambda \in \If_{47}$, then $\dim L(\lambda)=1$ and $E_i  v_\lambda= 0$, $F_i  v_\lambda= 0$, $g\gpl   v_\lambda= \lambda(g\gpl)v_\lambda$.
\end{lemma}
\pf
Let $N'(\lambda)=W(\lambda)+ W_1(\lambda)$. 
By a direct computation, $N'(\lambda)= \sum_{\beta\neq 0} M(\lambda)_\beta = N(\lambda)$. Therefore $L'(\lambda)=M(\lambda)/N'(\lambda)$ is
one-dimensional and simple.
\epf

\begin{example} Take $\Lambda=\mathbb Z_{12}=\langle g_2\rangle$, $g_1=g_2^8$ and $\gpl_1,\gpl_2\in\VLambda$ such that
	\begin{align}
	\gpl_1(g_2)&=\zeta^{11},& \gpl_2(g_2)&=-1;& &\text{hence} & \gpl_1(g_1)&=\zeta^4, &  \gpl_2(g_1)&=1.
	\end{align}
Applying the Main Theorem  we see that there is one simple module of dimension one and 
exactly $\#$ different isoclasses of a given dimension as in Table \ref{tab:quantity}:
\begin{table}[h]
	\caption{ Quantity of simple modules of dimension $>1$}\label{tab:quantity}
	\begin{center}
		\begin{tabular}{| c |c | c|c | c |c | c |c |} 
			\hline
			$\#$ & dimension & $\#$ & dimension & $\#$ & dimension & $\#$ & dimension
			\\ \hline
			67 & 144 &   7 & 108 &   10 & 96 & 		   2 & 85\\	\hline 
			6 & 72 &    4 & 71 & 		   4 & 61 &    2 & 49 \\	\hline
			10 & 48 &  4 & 47 &    6 & 37 &    7 & 36 \\	\hline
			4 & 35 &    4 & 25 &    2 & 23 & 4 & 11    \\	\hline
		\end{tabular}
	\end{center}
\end{table}

 Note that $\If_6$ and $\If_{10}$ are empty.
\end{example}

\end{document}